\documentclass[11pt,twoside]{article}
\usepackage{latexsym}
\usepackage{amssymb,amsbsy,amsmath,amsfonts,amssymb,amscd}
\usepackage{epsfig, graphicx}
\usepackage{color}
\usepackage{hyperref}
\newcommand  \ind[1]  {   {1\hspace{-1.2mm}{\rm I}}_{\{#1\} }    }

\newcommand{\commentout}[1]{}
\newcommand{\R}{\mathbb{R}}
\newcommand{\N}{\mathbb{N}}

\newcommand {\sg} {\sigma}

\newcommand {\Chi} {{\bf \raise 2pt \hbox{$\chi$}} }

\newcommand {\dv}  { {\rm div} }
\newcommand {\f}   {\frac}
\newcommand {\p}   {\partial}
\newcommand{\dis}{\displaystyle}

\newcommand{\beq}{\begin{equation}}
\newcommand{\eeq}{\end{equation}}
\newcommand{\bea} {\begin{array}{rl}}
\newcommand{\eea} {\end{array}}
\newcommand{\bepa}{\left\{ \begin{array}{l}}
\newcommand{\eepa} {\end{array}\right.}

\newcommand{\tildeu}{\widetilde{u}}
\newcommand{\tildec}{\widetilde{c}}
\newtheorem{theorem}{Theorem}[section]
\newtheorem{lemma}[theorem]{Lemma}

\newtheorem{remark}[theorem]{Remark}
\newtheorem{proposition}[theorem]{Proposition}

\newcommand{\qed}{{ \hfill
                       {\unskip\kern 6pt\penalty 500 \raise -2pt\hbox{\vrule\vbox to 6pt{\hrule width 6pt
                       \vfill\hrule}\vrule} \par}   }}

\title{\Large \bf Traveling wave solution of the Hele-Shaw model of tumor growth with nutrient}

\author{
Beno\^it Perthame\thanks{Sorbonne Universit\'es, UPMC Univ Paris 06, UMR 7598, Laboratoire Jacques-Louis Lions, F-75005, Paris, France} \thanks{CNRS, UMR 7598, Laboratoire Jacques-Louis Lions, F-75005, Paris, France}
\thanks{INRIA-Paris-Rocquencourt, EPC MAMBA, Domaine de Voluceau, BP105, 78153 Le Chesnay Cedex, France} \footnotemark[5]
\and
Min Tang\thanks{Department of mathematics and Institute of Natural Sciences. MOE-LSC, Shanghai Jiao Tong University, China} \footnotemark[5]
\and 
Nicolas Vauchelet\footnotemark[1] \footnotemark[2] \footnotemark[3]
\thanks{email~: benoit.perthame@upmc.fr, tangmin@sjtu.edu.cn, nicolas.vauchelet@upmc.fr}
}

\date{\today}

\begin{document}
\maketitle
\pagestyle{plain}
\begin{abstract}

Several mathematical models of tumor growth are now commonly used to explain medical observations and predict cancer evolution based on images. These models incorporate mechanical laws  for tissue compression combined with rules for nutrients availability which can differ depending on the situation under consideration, {\it in vivo } or {\it in vitro}. Numerical solutions exhibit, as expected from medical observations, a proliferative rim  and a necrotic core. However, their precise profiles are rather complex, both in one and two dimensions. 

We study a simple  free boundary model formed of  a Hele-Shaw equation for the cell number density  coupled to a diffusion equation for a nutrient. We can prove that a traveling wave solution exists with a healthy region separated from the progressing  tumor by a sharp front (the free boundary) while the transition to the necrotic core is smoother.  Remarkable is the pressure distribution which vanishes at the boundary of the proliferative rim with a vanishing derivative at the transition point to the necrotic core.

\end{abstract}

\\[5pt]
\noindent {\bf Mathematical Classification numbers:} 76D27; 35K57;  35C07; 92C50;
\pagenumbering{arabic}

\section{Introduction}

Mathematical models of tumor growth, based on a mechanistic approach, have been developped and studied in many works. Nowadays, they are being used for image analysis and medical predictions \cite{Swanson, bresch,cornelis}.  Among these models, we can distinguish two main directions. Either the dynamics
of the cell population density is described at the cell level 
using fluid mechanical concepts, or a geometrical description is used leading to a free boundary problem 
\cite{bellomo1, bellomo2, friedman, friedman_hu,greenspan,Lowengrub_survey, chapman_maini}.
The link between these two approaches has been established in \cite{PQV}, using the asymptotic of a stiff law-of-state for the pressure and it has been extended to the case with active cells in \cite{PQTV}.
\\

A simple cell population density model governing the time dynamics
of the cell population density $n(x,t)$ under pressure forces and cell multiplication governed by the local availability of nutrients with concentration $c(x,t)$, writes~:
\beq
\f{\p}{\p t} n - \dv \big( n \nabla p(n) \big) = n G\big(c\big),  \qquad x\in \R^d,\; t \geq 0,
\label{eq:pqn}
\eeq
where the nutrient dependent term $G$ represents  growth. It can be used with $d=2$ or $3$, for laboratory experiments or  {\it in vivo} representation. The pressure law is given by~: 
$$
p(n) = n^\gamma.
$$
We assume that the nutrient is single and that diffusion and consumption can be described by an elliptic  equation (they are fast compared to the time scale of cell division)~:
\beq
- \Delta c + \Psi(n,c) =0.
\label{eq:genc}
\eeq
Here the function $\Psi(n,c)$ describes both the effect of the vasculature network bringing the 
nutrient to the cells and nutrient consumption by the cells. At this stage, the reaction terms share similarities with a classical  two-component chemical reaction system for reactant and temperature \cite{BNS, BBH, Logak} and also with models of bacterial swarming \cite{BenJ, mimura}.
\\

In this paper, we consider two models for the nutrient consumption.
For the  {\it in vitro} model, we assume that the tumor is surrounded by
a liquid in which the nutrient diffusion is so fast that
it is assumed to be constant; while inside the tumoral region, 
the consumption is linear, i.e. $\Psi(n,c)=\psi(n)c$, with $\psi$ a 
nonnegative function. Thus, for {\it in vitro} models,  the equation for the nutrient consumption
writes
\beq\label{eq:vivoc}
-\Delta c+ \psi(n)c=0, \quad \mbox{ for } x\in \{n>0\}\ ;
 \qquad c=c_B, \quad \mbox{ for } x\in \{n=0\}.
\eeq
For the model {\it in vivo}, we consider that the nutrient is brought by the vasculature 
network in the healthy region $\{n=0\}$ and diffused to the tissue. We choose
$\Psi(n,c)=\psi(n)c+\ind{n=0}(c-c_B)$ and the system writes 
\beq\label{eq:vitroc}
-\Delta c + \psi(n) c = \ind{n=0}(c_B-c).
\eeq
Here $\ind{n=0}(c_B-c)$ indicates that the vasculature is  pushed away by the growing tumor so that
the nutrient is directly available only from healthy tissues; a full discussion of this issue can be found in \cite{CBCB,nick_phd}.

The asymptotic limit $\gamma\to +\infty$ can be viewed as an incompressible
limit. It leads to the free boundary problem  of Hele-Shaw type. Since when $n<1$, we have $p(n) \underset{\gamma \to \infty}{\longrightarrow} 0$, we can write the Hele-Shaw equation in a {\em weak form} as  
\beq \label{as:tgnnlimit}
\left\{\begin{array}{l}
\dis \f{\p}{\p t} n - \dv \big( n \nabla p \big) = n G\big(c\big),  
\qquad x\in \R^d,\; t \geq 0, 
\\[2mm]
\dis n=1, \quad \text{ for }  x\in \Omega(t) = \{p(t)>0\}, 
\end{array} \right.
\eeq
\beq \label{as:tgnclimit}
 - \Delta c + \Psi(n,c) =0,\qquad x\in \R^d.
\eeq
But one can also establish a {\em strong form} as a free boundary (geometrical) problem. Multiplying \eqref{eq:pqn} by $p'(n)$, we obtain the  following equation on $p$~:
$$
\f{\p}{\p t} p- n p'(n) \Delta p + |\nabla p|^2 = n p'(n)  G\big(c\big)
$$
and for the special case $p= n^\gamma$ at hand, we find
$$
\f{\p}{\p t} p- \gamma p \Delta p + |\nabla p|^2 = \gamma p G\big(c\big).
$$
Letting formally $\gamma \to \infty$ we find 
$$
- p \Delta p = p G\big(c\big).
$$
The Hele-Shaw geometrical model (see \cite{PQV} for more details) is to write 
\beq \label{as:tgnplimit}
\left\{\begin{array}{rll}
-\Delta p &=  G(c) \qquad &\hbox{in }\, \Omega (t),
\\[2mm]
p&=0 \qquad &\hbox{on }\p \Omega (t), 
\end{array}\right.
\eeq
together with the velocity of the free boundary $\partial \Omega(t)$, $v =  - \nabla p$.
\\

Equation \eqref{as:tgnclimit} describes consumption and diffusion of the
nutrient through tumoral tissue. 
In equation \eqref{as:tgnnlimit}, cells multiplication is limited by nutrients
brought by blood vessels. This dependency is described 
by the increasing function $G(\cdot)$.
Moreover, lack of nutrients leading to cells necrosis is modeled
by assuming that $G$ can take negative values~:
\beq \label{as:tgnbis}
G'(\cdot) >0, \qquad  G(\bar c) = 0 \quad \text{ for some } \bar  c >0 .
\eeq

This paper is devoted to the description of the structure of the solutions
of the 'incompressible' model, that is the Hele-Sahw equation. This question arises because numerical observations, which we present in Section~\ref{sec:num},  show that the solutions are rather complex, exhibiting a sharp front and an apparently smooth transition to a necrotic core. This type of shape  is comparable to agent based simulations and experimental observations in \cite{ByDr, Nick}.
In Section~\ref{sec:TW}, we investigate the existence of traveling waves in one dimension for the case {\it in vitro}. We establish in particular that 
the waves are formed by a proliferative zone, where the density of tumoral 
cells is maximal, followed by a necrotic zone where the cell density decreases 
towards zero. This construction is extended to the case {\it in vivo} in Section~\ref{sec:invivo}.

\section{Numerical observations}
\label{sec:num}

In two dimensions, numerical simulations of system \eqref{eq:pqn},  \eqref{eq:genc} with $\gamma$ large, exhibit various types of patterns that we present now and that motivate to conduct a detailed analysis. For instance, a more accurate view of the behaviour of solutions can be obtained with one dimensional solutions and they still exhibit a complex structure.

\subsection{Two dimensional simulations}

For two dimensional simulations,  we have considered the equation on nutrient given by  
$$
-\Delta c+ \psi(n)c=0, \qquad (x,y)\in\Omega.
$$ 
The parameters and nonlinearities are chosen as 
$$
\gamma=30,\qquad \psi(n)=4 n ,\qquad G(c)=200(c-0.3),
$$
the initial data is given by
$$
n(x,y)=\left\{\begin{array}{ll}
0.7,& \qquad 0 \leq \sqrt{x^2+y^2}<0.8,
\\[5pt]
 0.7*\exp\big(-20(x^2+y^2-0.64)\big),& \quad \; 0.8<\sqrt{x^2+y^2}<2,
\\[5pt]
0,& \qquad 2<\sqrt{x^2+y^2}<3,
\end{array}\right.
$$

The numerics uses the finite element method implemented within the software FreeFem++  \cite{FreeFem}. The elliptic equation for $c$ is discretized thanks to P1 finite element method.
For $n$, we use a time splitting method by first solving
$$
\f{\p}{\p t} n - \dv \big( n \nabla p(n) \big) = 0
$$
with P1 finite element method for one time step; 
and then solving
$$
\f{\p}{\p t} n = n G\big(c(x,t) \big),
$$
for one time step to update the value of $n$ for the next time step. 
The computational domain is a disc with radius $3$ and we denote by $\Gamma$ its
boundary. The numerics  is obtained with the Dirichlet boundary conditions for both $c$ and $n$~:
$$
 c|_{\Gamma}=1,\qquad n|_{\Gamma}=0.
$$ 
This boundary condition on $c$ seems necessary for the instabilities we present below (see further comments below). With the models \eqref{eq:vivoc}, \eqref{eq:vitroc} we obtain radial symmetric solutions which we do not present.

\medskip

\begin{figure}
\begin{center}
\includegraphics[width=0.3\textwidth]{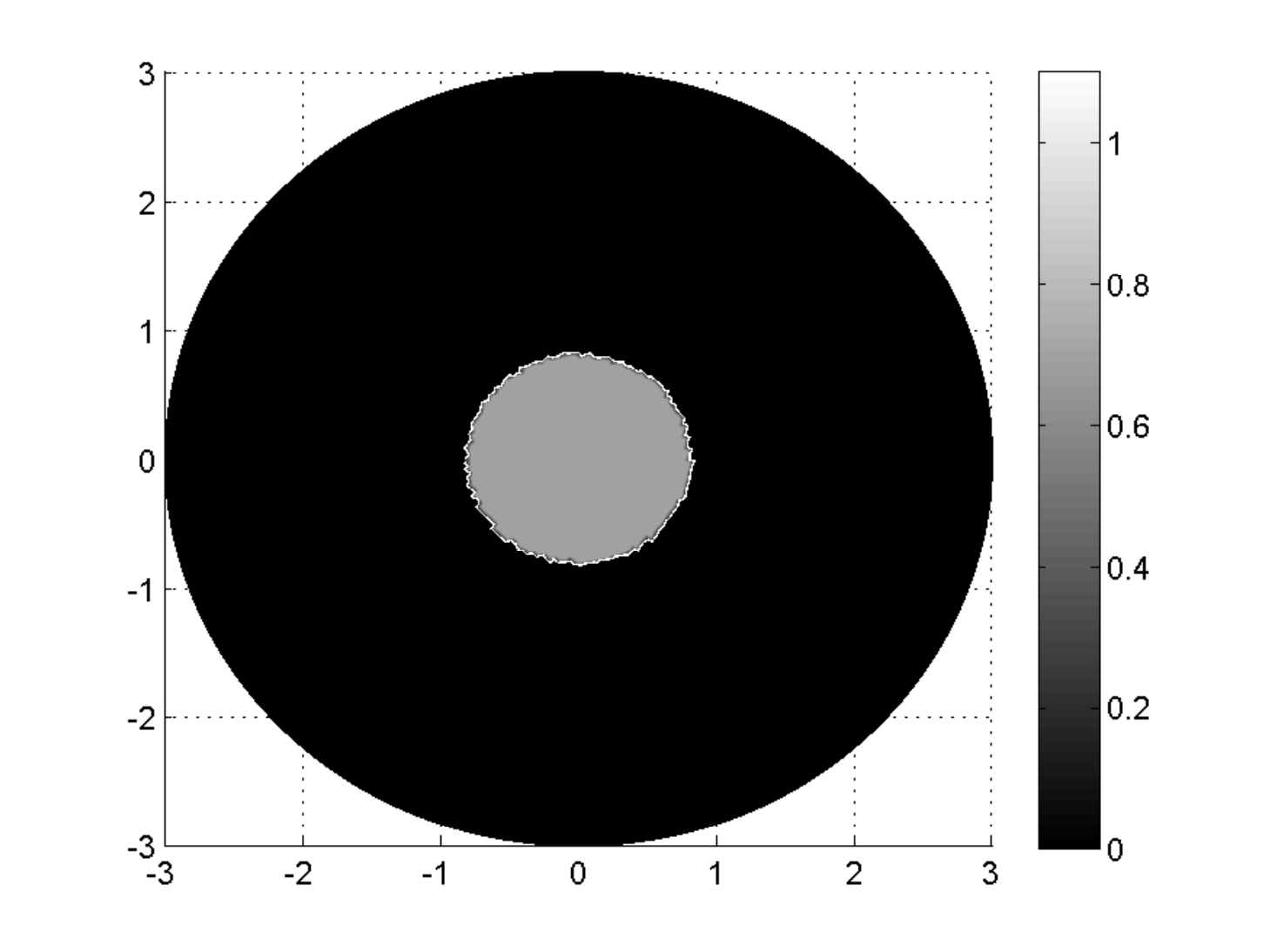}
\includegraphics[width=0.3\textwidth]{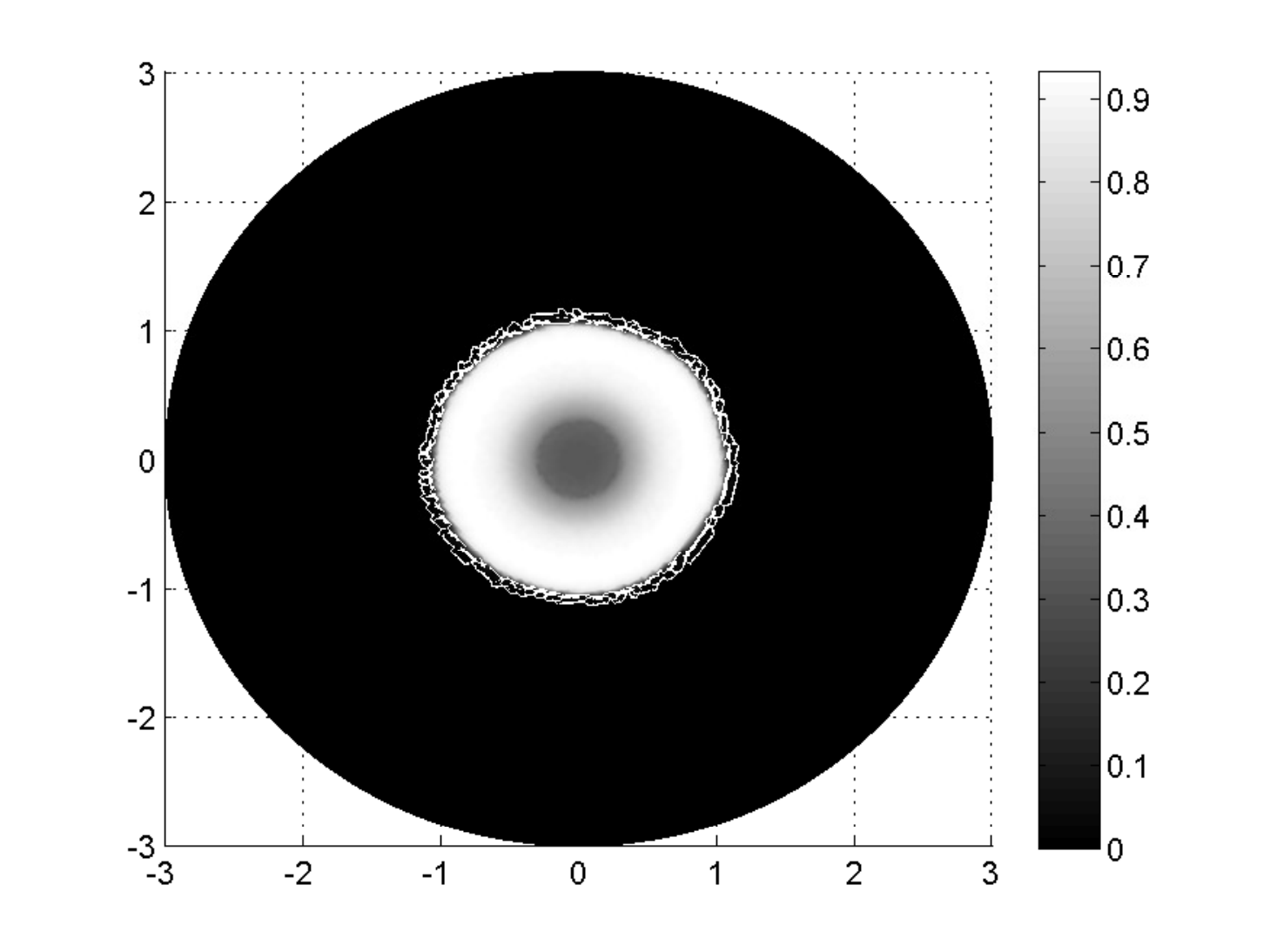}
\includegraphics[width=0.3\textwidth]{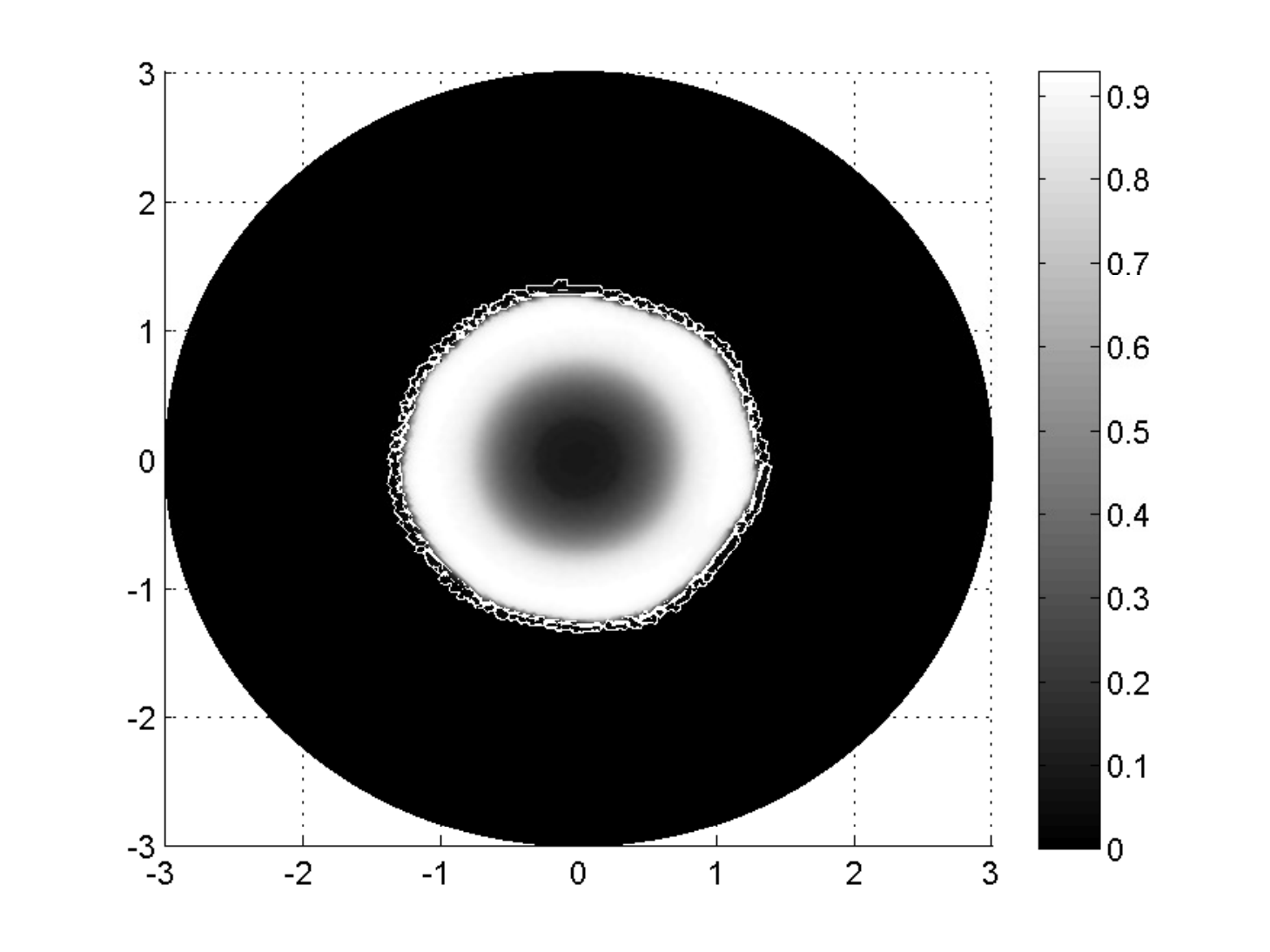}
\includegraphics[width=0.3\textwidth]{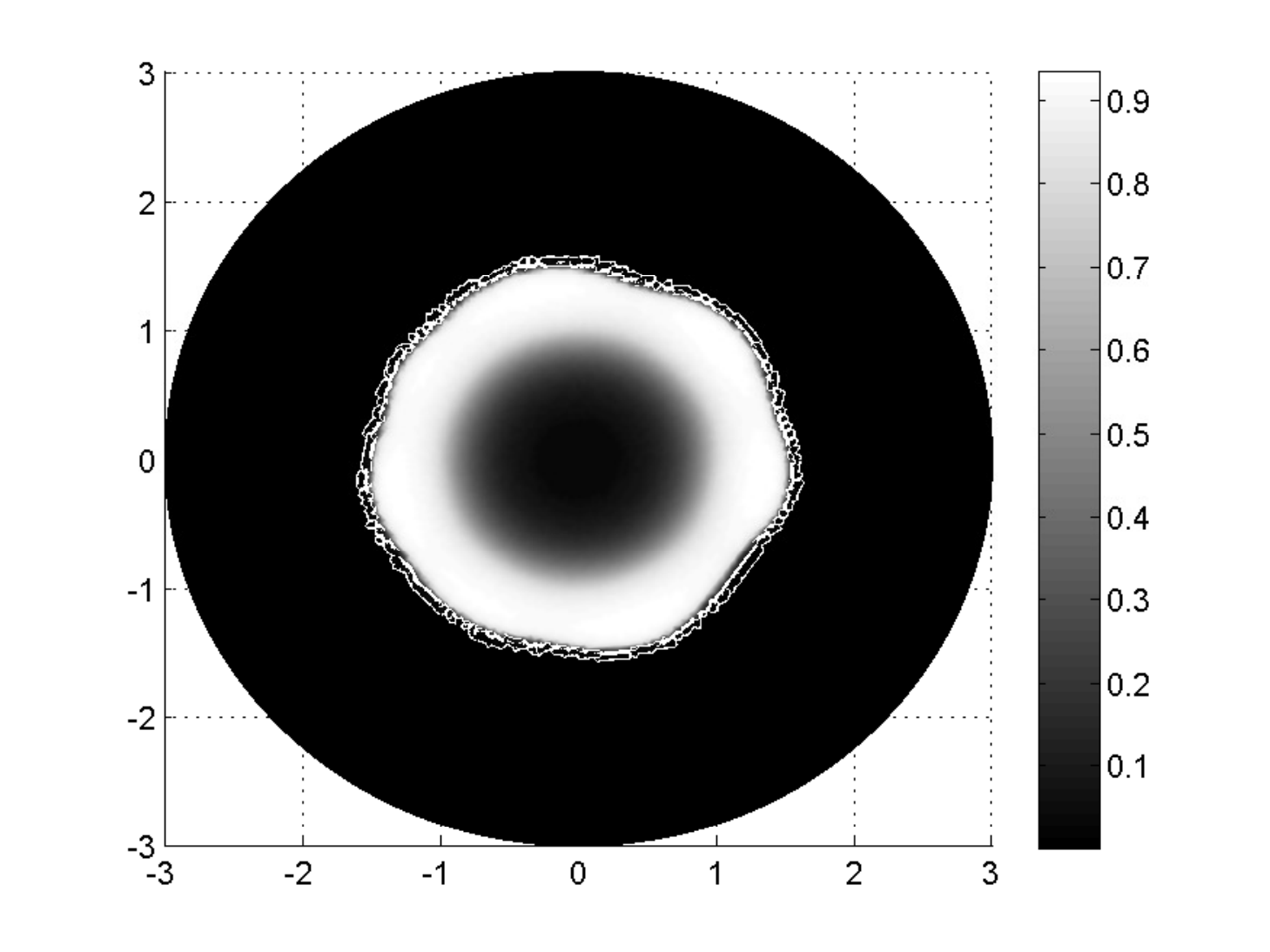}
\includegraphics[width=0.3\textwidth]{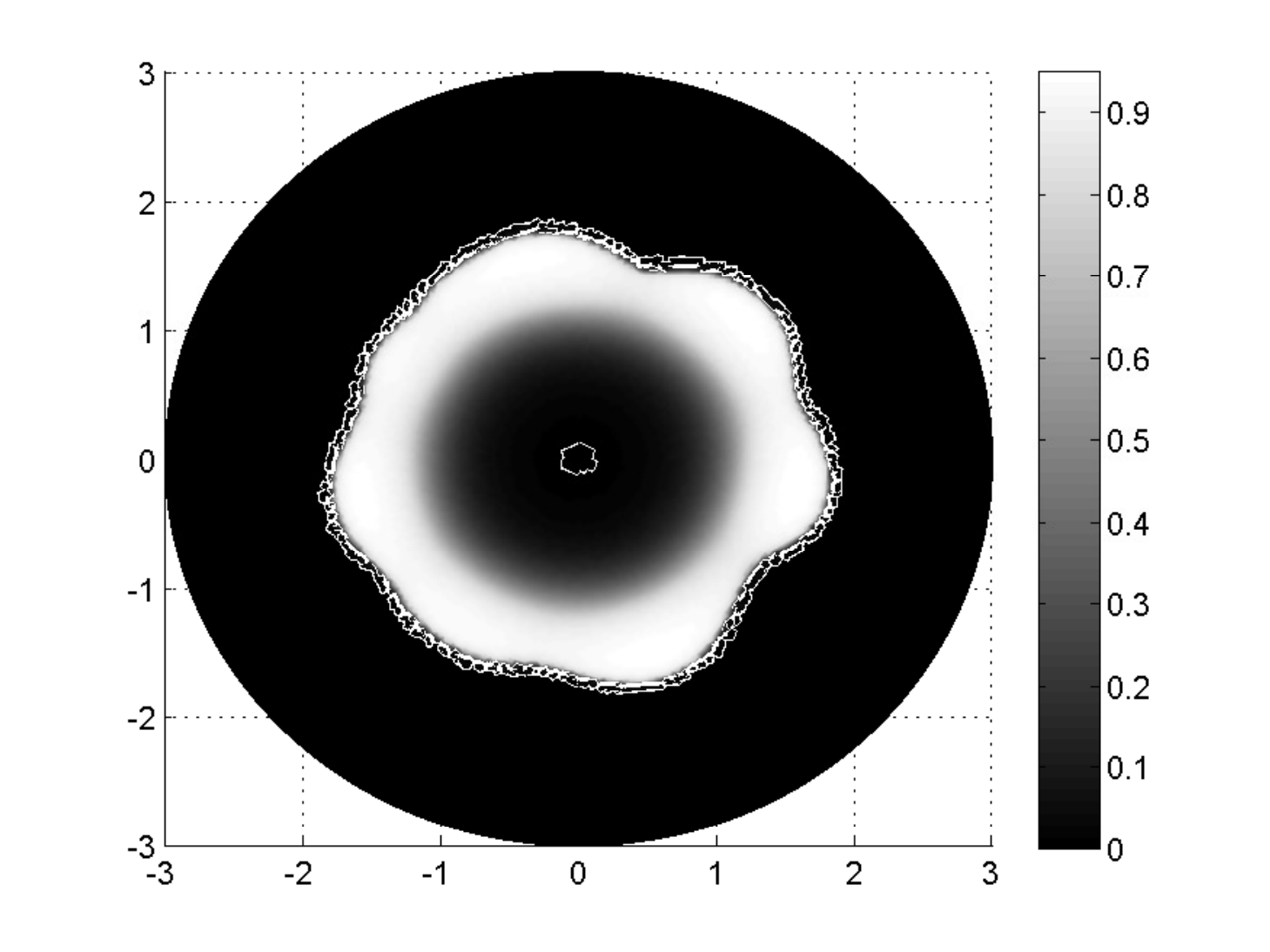}
\includegraphics[width=0.3\textwidth]{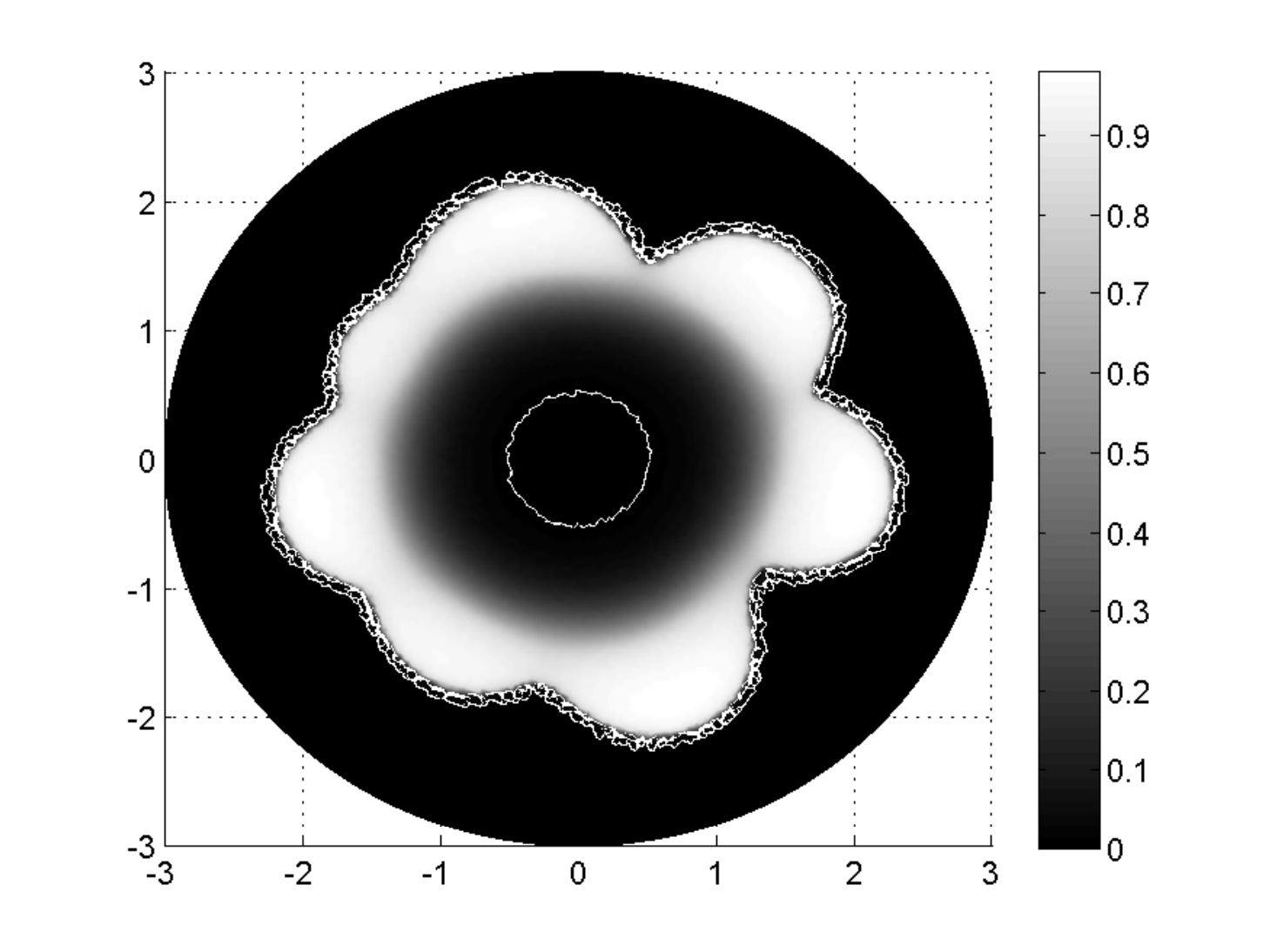}
\caption{The time dynamics of the density $n$ from initially radially symmetric plateau.
The time increases from left to right and from top to bottom. }
\label{fig:density}
\end{center}
\end{figure}

\begin{figure}
\begin{center}
\includegraphics[width=0.3\textwidth]{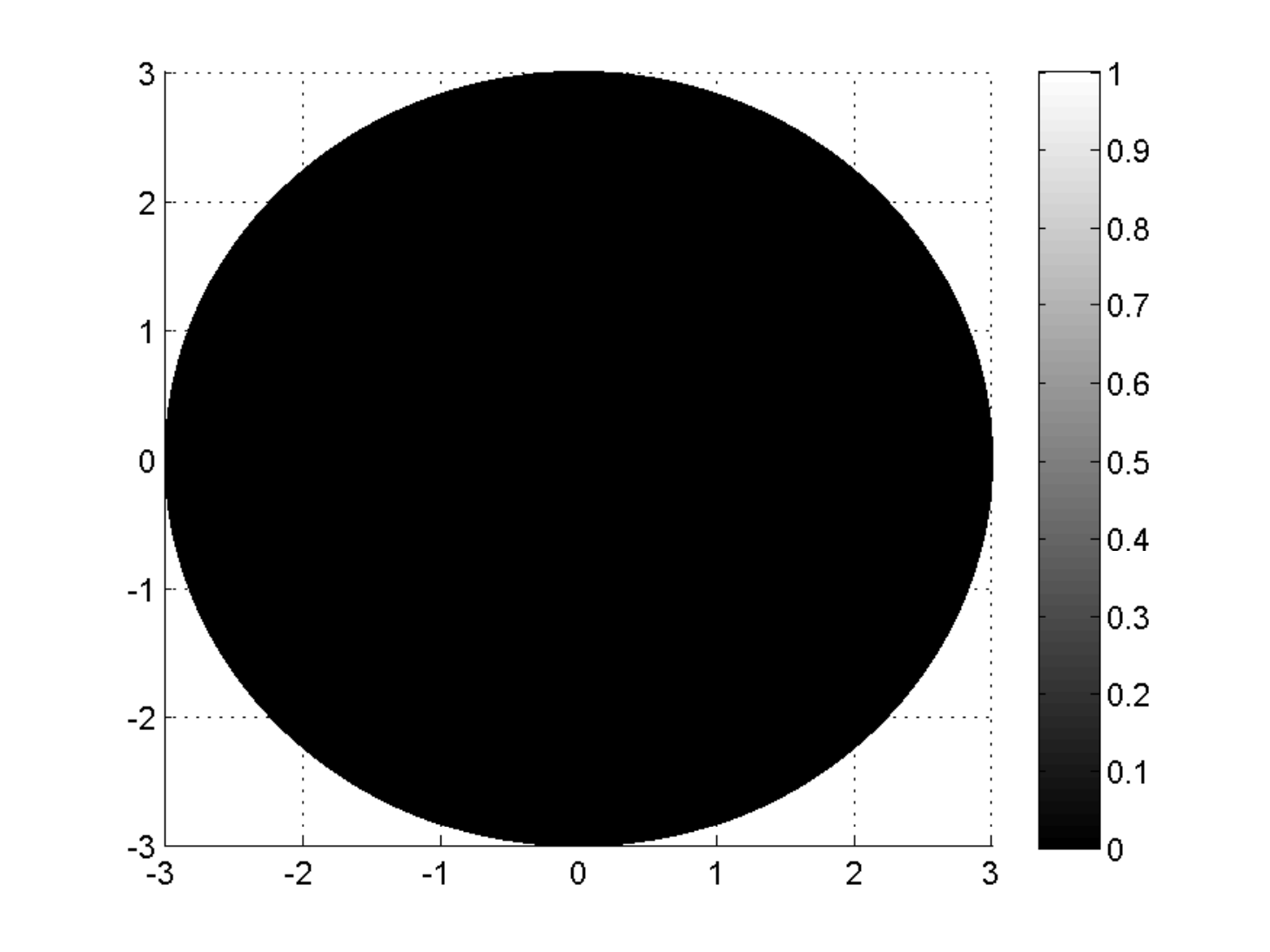}
\includegraphics[width=0.3\textwidth]{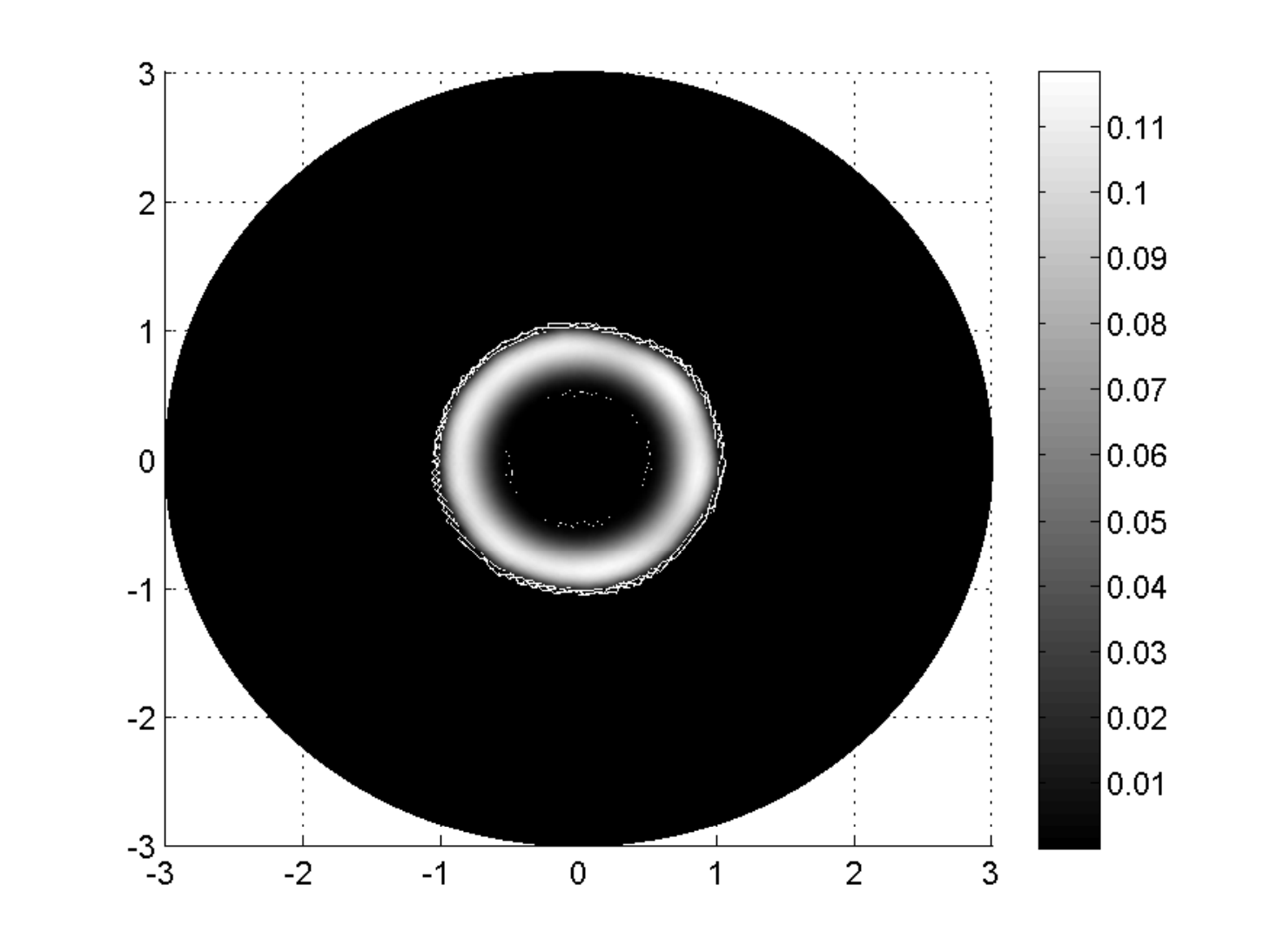}
\includegraphics[width=0.3\textwidth]{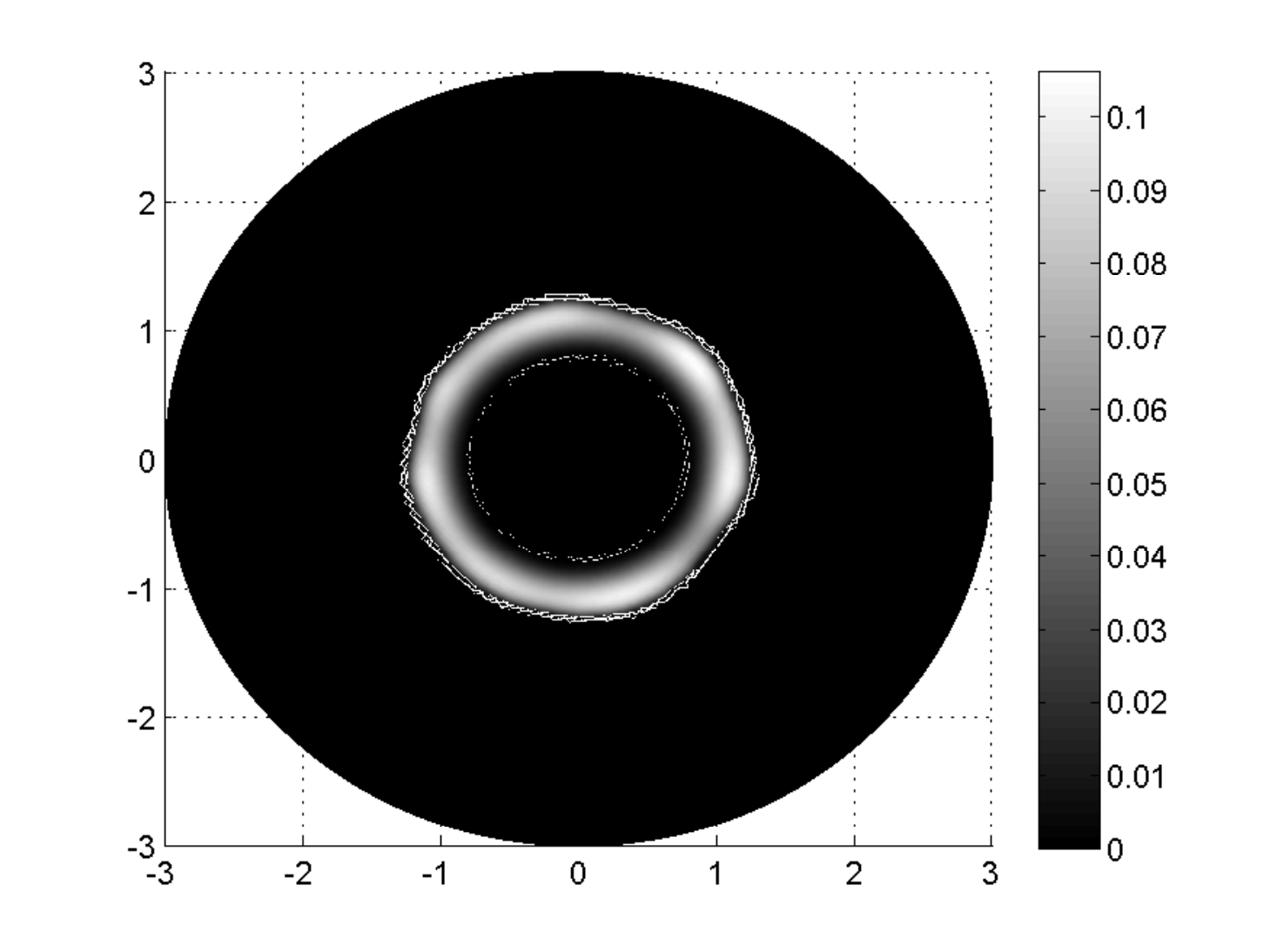}
\includegraphics[width=0.3\textwidth]{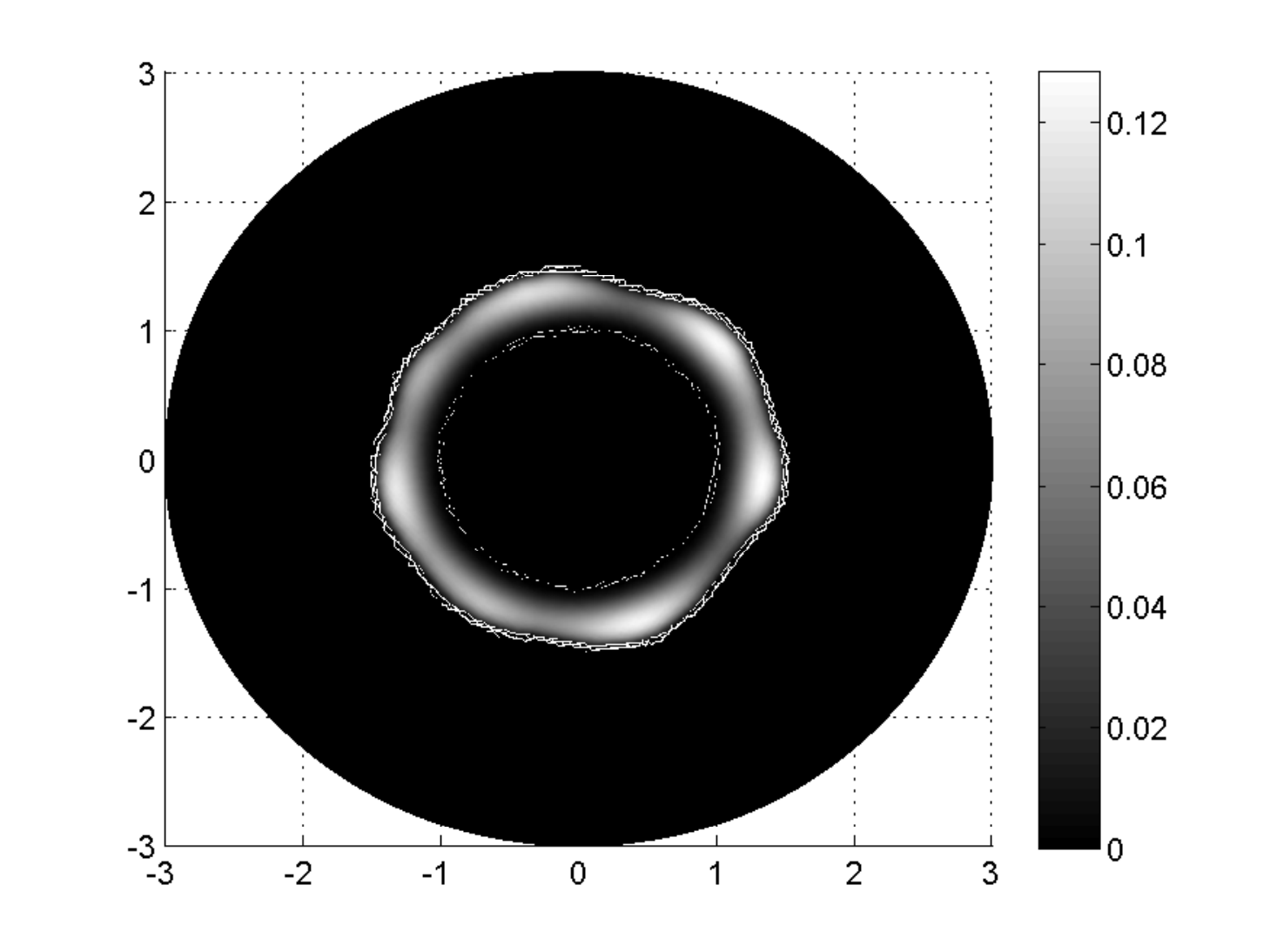}
\includegraphics[width=0.3\textwidth]{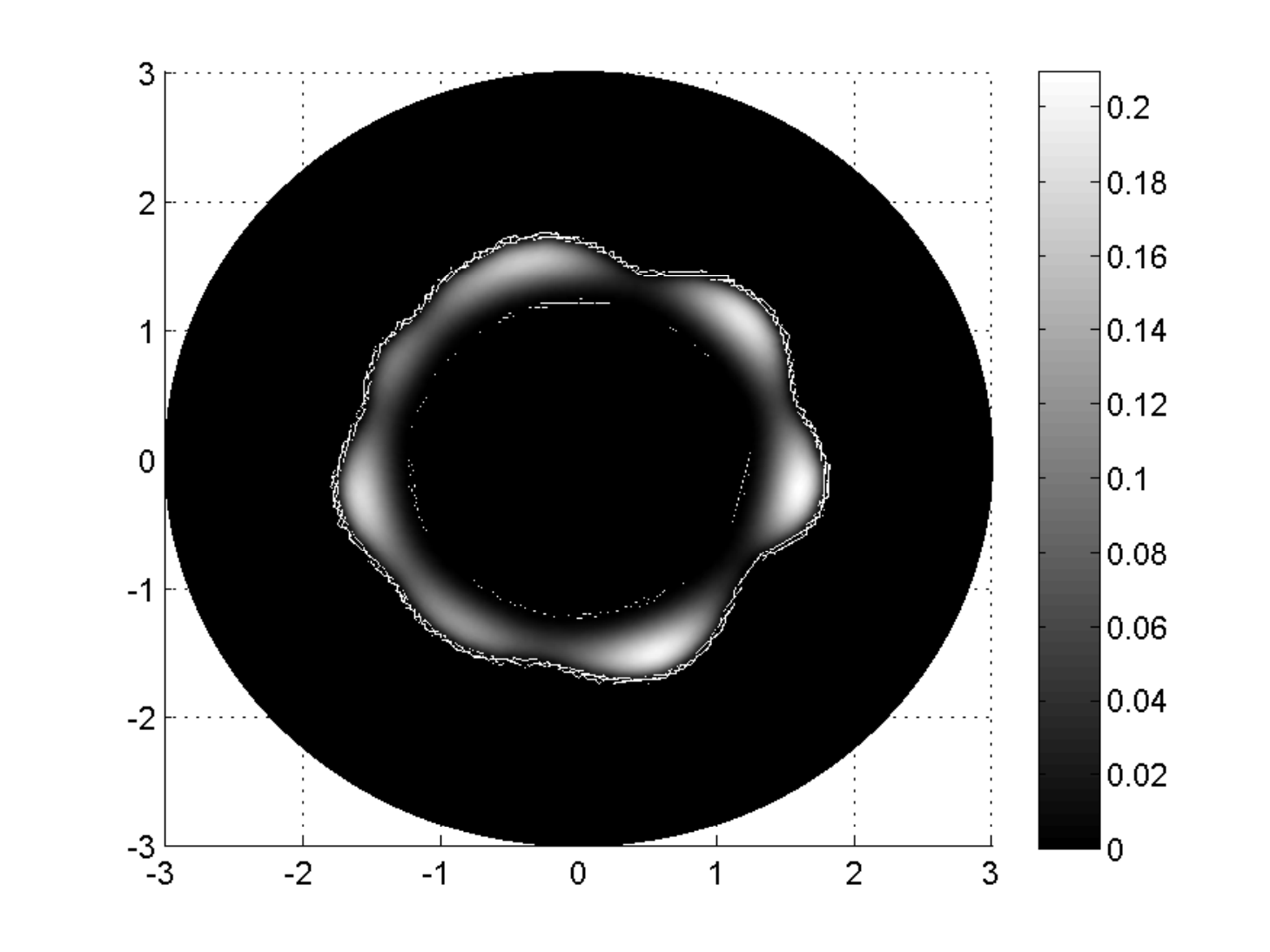}
\includegraphics[width=0.3\textwidth]{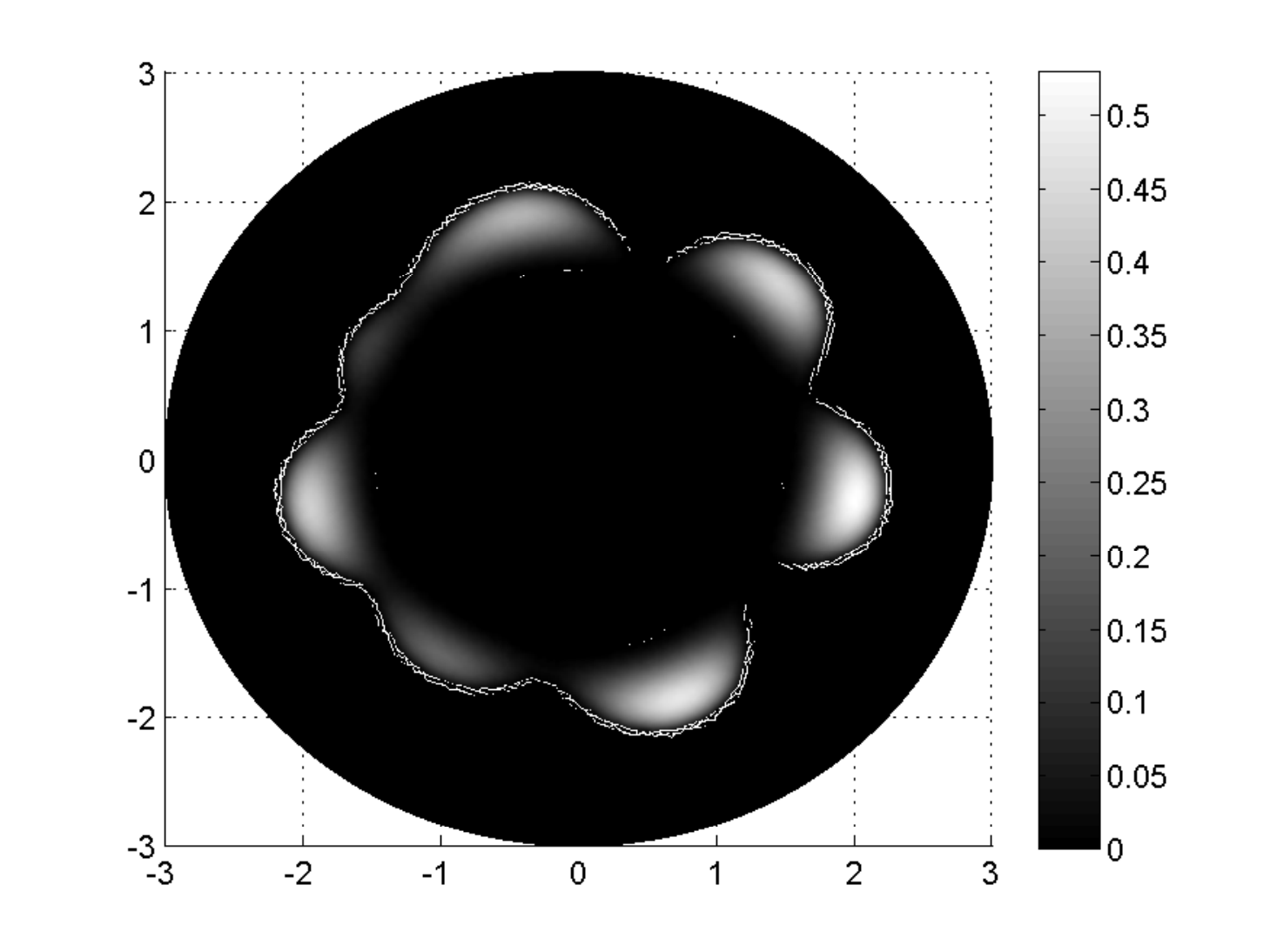}
\caption{The time evolution of the pressure $P(n)$ from initially radially symmetric plateau.
The time increases from left to right and from top to bottom. }
\label{fig:pressure}
\end{center}
\end{figure}

\begin{figure}
\begin{center}
\includegraphics[width=0.3\textwidth]{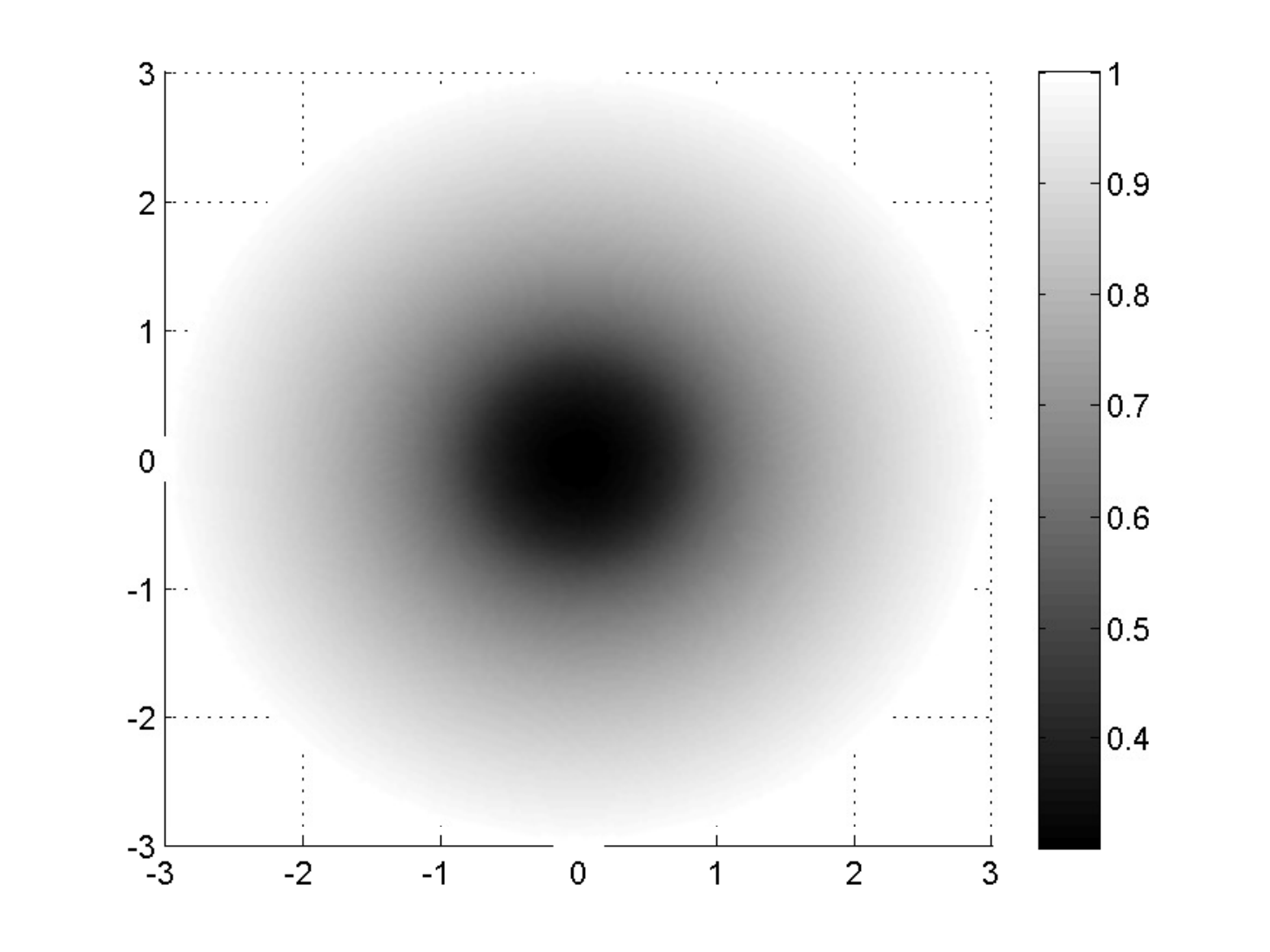}
\includegraphics[width=0.3\textwidth]{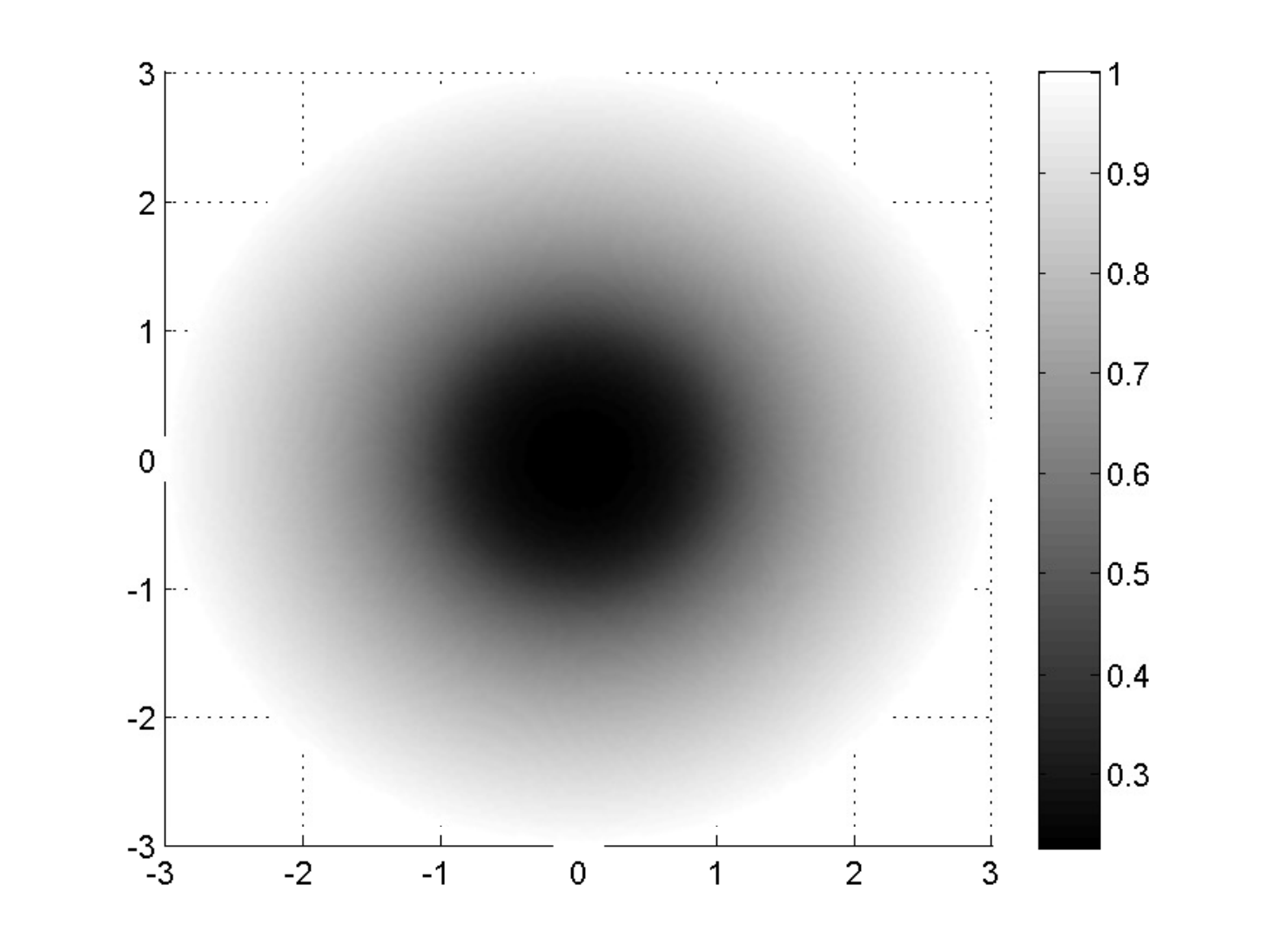}
\includegraphics[width=0.3\textwidth]{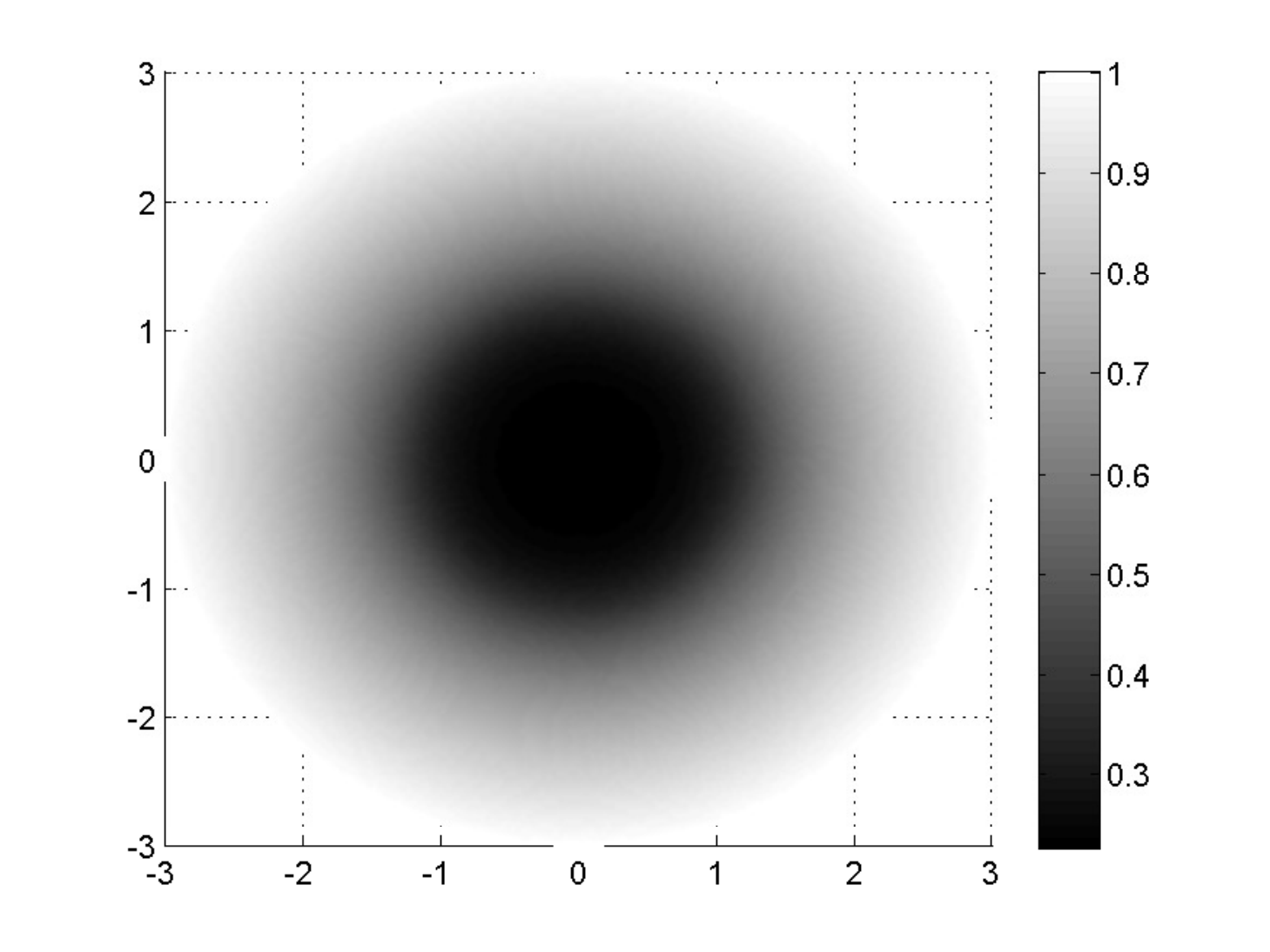}
\includegraphics[width=0.3\textwidth]{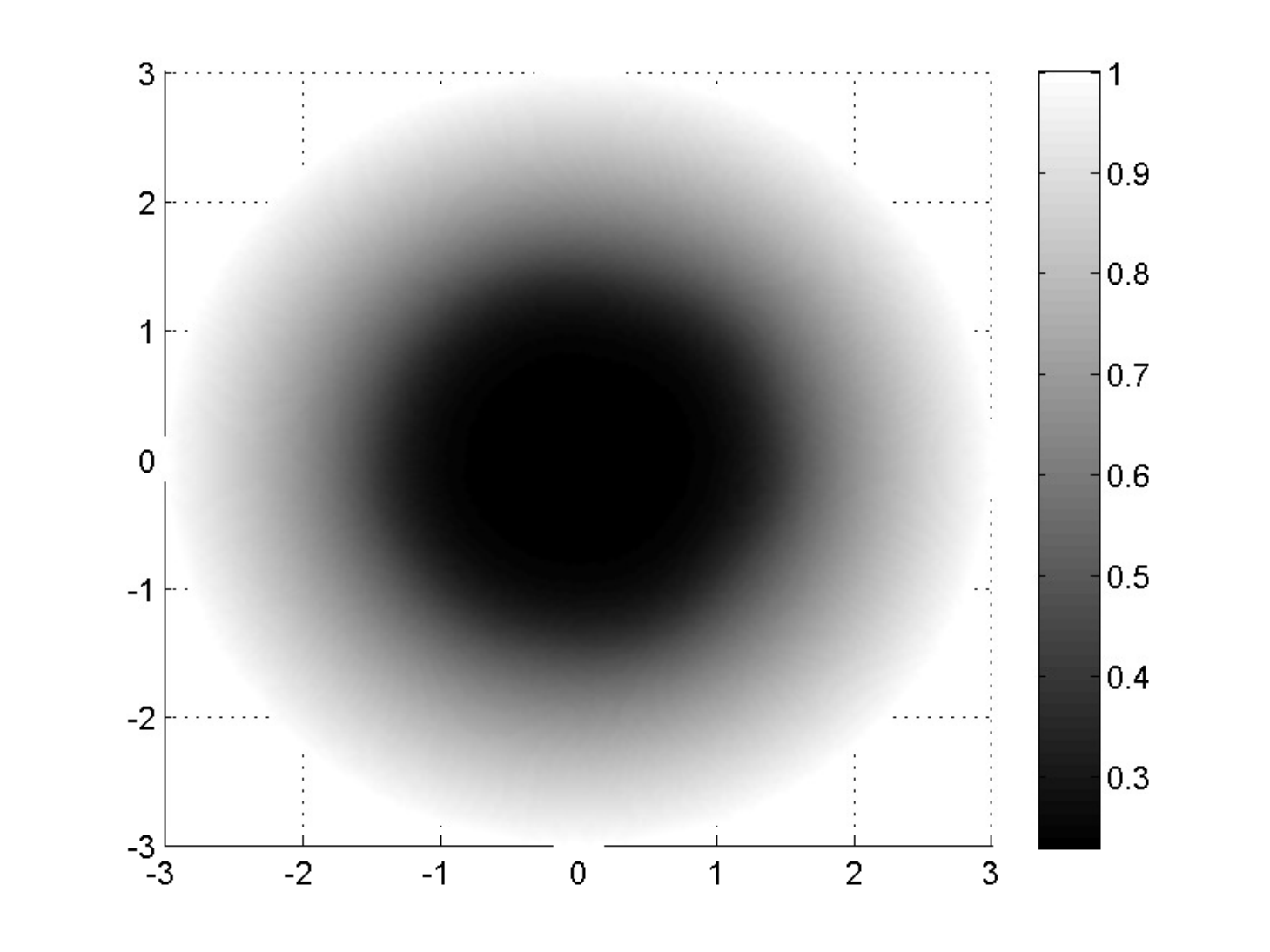}
\includegraphics[width=0.3\textwidth]{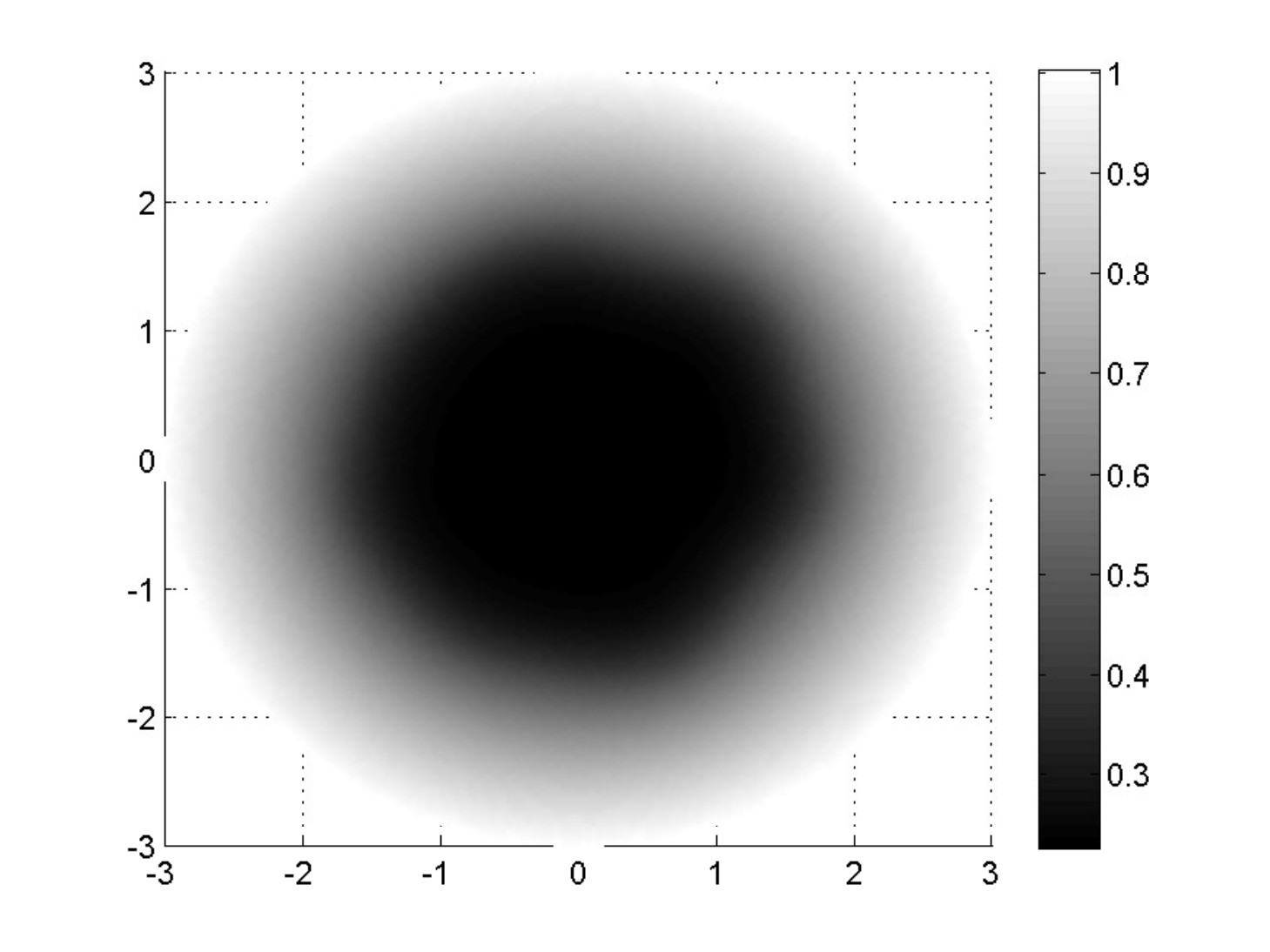}
\includegraphics[width=0.3\textwidth]{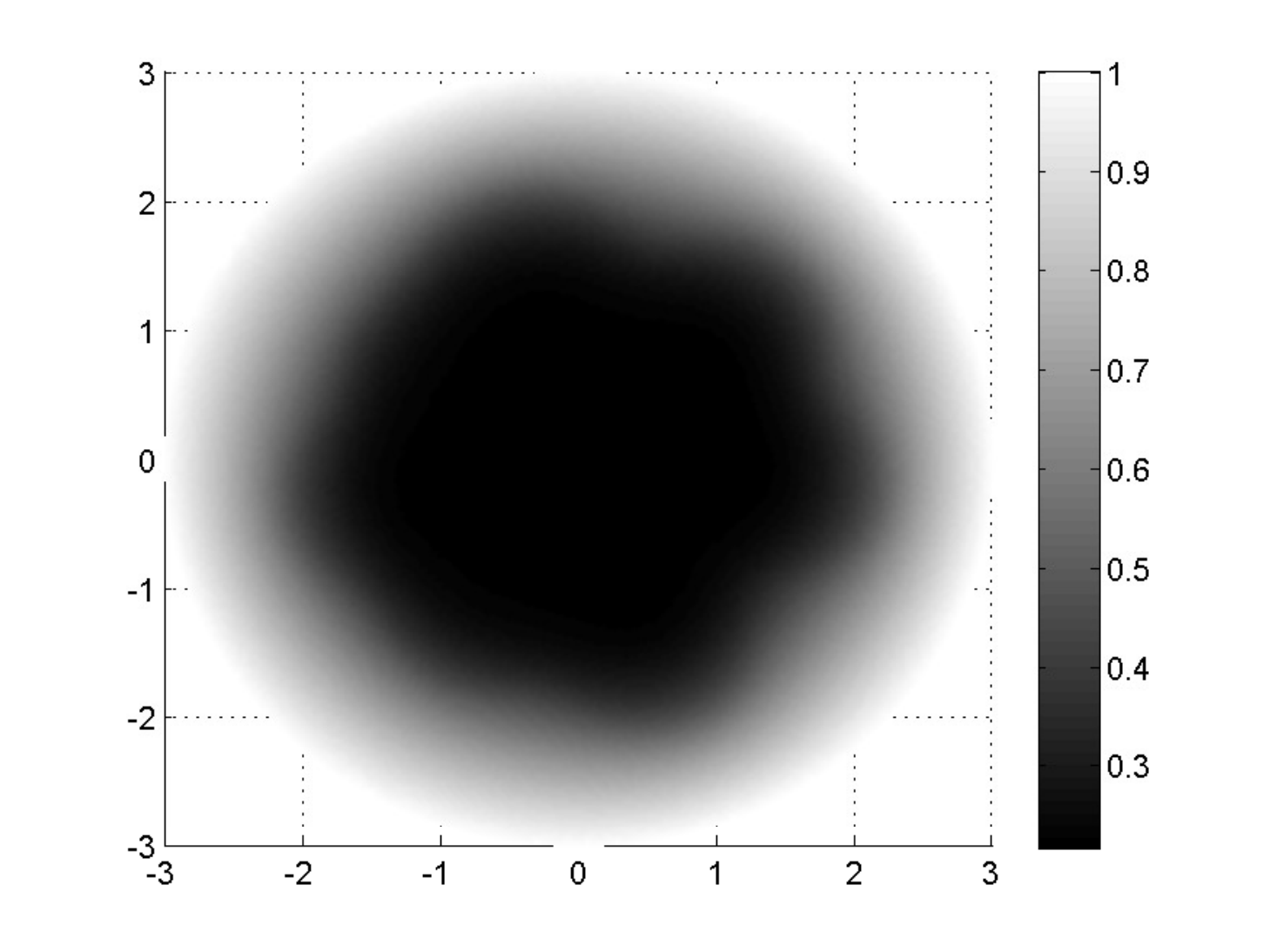}
\caption{The time evolution of the nutrient distribution $c$ from initially radially symmetric plateau.
The time increases from left to right and from top to bottom. }
\label{fig:nutrient}
\end{center}
\end{figure}

Fig.~\ref{fig:density},  Fig.~\ref{fig:pressure} and Fig.~\ref{fig:nutrient}
display the time dynamics for, respectively, the cell number density $n$, the pressure $p(n)$ and the nurient concentration $c$ .
\\

We can distinguish the different phases of the dynamics. We first observe 
a first step of growth where the cell number density $n$ builds up without motion. After this transitory step, the tumor invades the tissue with spherical symmetry. Finally,  a necrotic core can be formed at the center of the tumor
 where cells die due to the low level of nutrient, this is observed in Fig.~\ref{fig:density}. Then, the proliferating cell density is higher at the outer rim and decreases in the middle which corresponds to the inner necrotic core. In this phase, we can also observe symmetry breaking  
after some time. This is  related to spheroid instability in the tumor already observed by several authors \cite{BenAmar} and by \cite{Kessler_Levine} for a simpler, but related, system of two-component reaction-diffusion.


\subsection{One dimensional numerical traveling waves}

Because we wish to study the invasion process before the instability occurs, we focus on dimension $1$.
We present some numerical simulations for the case  {\it in vitro} given by equation \eqref{eq:vivoc}  with the following parameters
\beq\label{eq:ex1parameter}
\gamma=50,
\qquad \psi(n)=\left\{\begin{array}{ll}
2,&n\geq1,\\
1,&n<1,
\end{array}\right.\qquad
G(c)=\left\{\begin{array}{ll} \ \ 21,&c\geq 0.6,\\
-30,&c<0.6.\end{array}\right.
\eeq
The computational domain is $[-5,5]$ and we start from an initial plateau 
$$
n(x)=\left\{\begin{array}{ll}0.1, & x\in(-0.5,0.5), 
\\
\ \ 0, & x\in [-5,-0.5]\cup[0.5,5].
\end{array}\right.
$$
The boundary condition are as follows:
$$
n(-5)=n(5)=0,\qquad c(-5)=c(5)=c_B.
$$

In order to track the front position, that is the boundary point of the set $\{n>0\}$,
and to perform the numerical simulations on the whole domain, we consider different diffusion values inside and outside the tumor region. 
Therefore, we assume that the diffusion coefficient for $c$ 
depends on $n$ in such a way that, at time $t$, the nutrient distribution satisfies
\beq\label{eq:phic}
-\nabla(H(n,t)\nabla c)+\phi(n,t)c=\phi(n,t)(1-H(n,t)) c_B.
\eeq
Here $H(n,t)$ gives the front position of $\{n(t)>0\}$  and 
$$
\phi(n,t)=\left\{\begin{array}{ll}
\psi(n),& \text{for } \; H(n,t)=1,
\\
1,&  \text{for } \; H(n,t)=0.\end{array}\right.
$$ 
More precisely, let $\epsilon$ be some small
tolerance (we choose $10^{-5}$ in practice), numerically $H(n,t)$ can be approximated by
$$
H(n,t)=\left\{\begin{array}{ll}1,&\quad\{n(s)>\epsilon, \mbox{ for some $s \in[0,t]$}\},\\
0,&\quad\mbox{otherwise}.\end{array}\right.
$$
Therefore, the nutrient diffusion coefficient is $1$ inside the tumor while it is worth $0$ outside. Formally, since $\psi(n)\neq 0$ according to \eqref{eq:ex1parameter}, the equation 
\eqref{eq:phic}
is equivalent to write
$$
\left\{\begin{array}{ll}-\Delta c+\lambda\Psi(n)c=0,&\quad  \mbox{for } H(n,T)=1,
\\
c=c_B, & \quad  \mbox{for } H(n,T)=0.
\end{array}\right.
$$

\begin{figure}
\begin{center}
a)\includegraphics[width=0.45\textwidth]{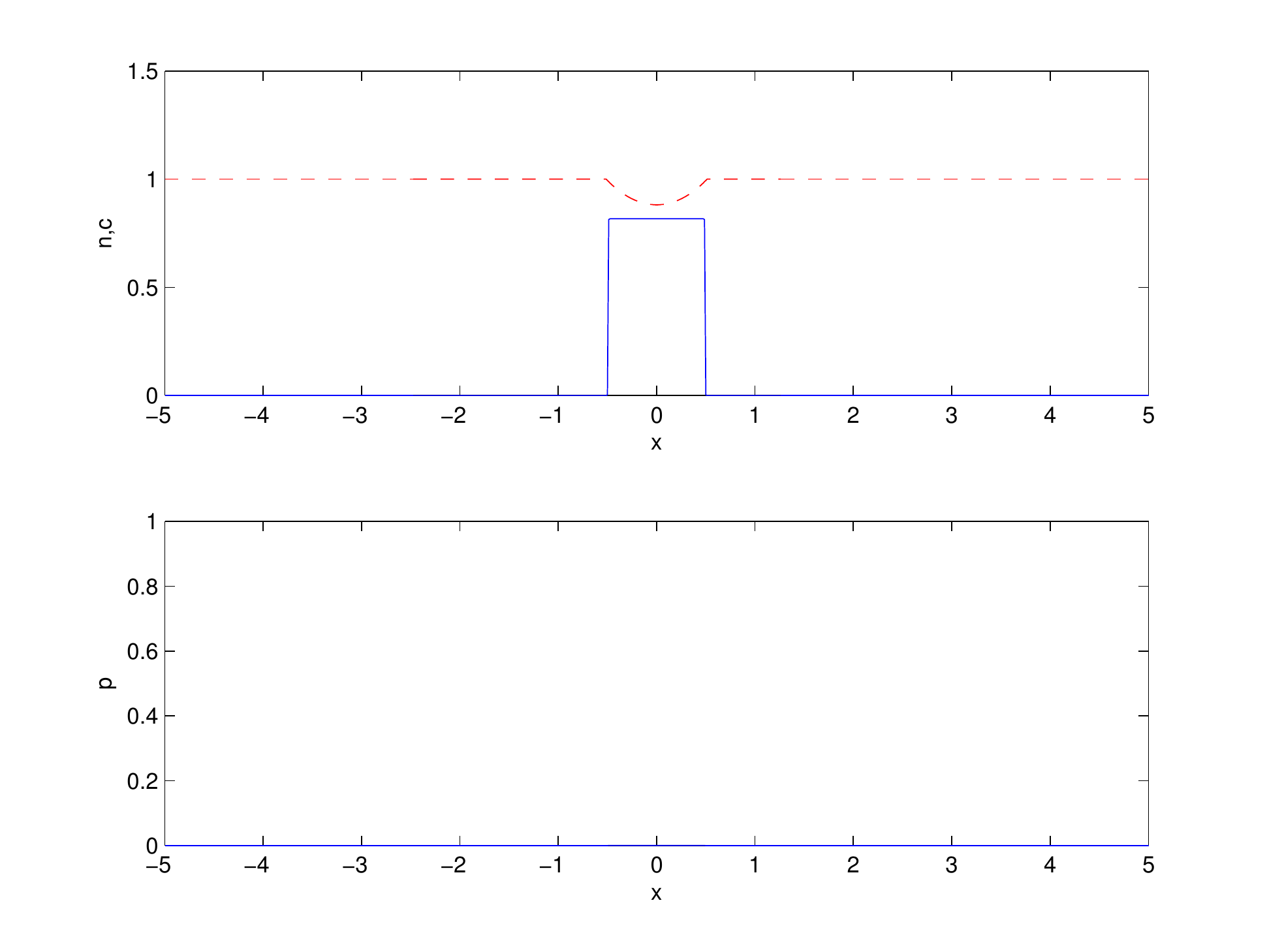}
b)\includegraphics[width=0.45\textwidth]{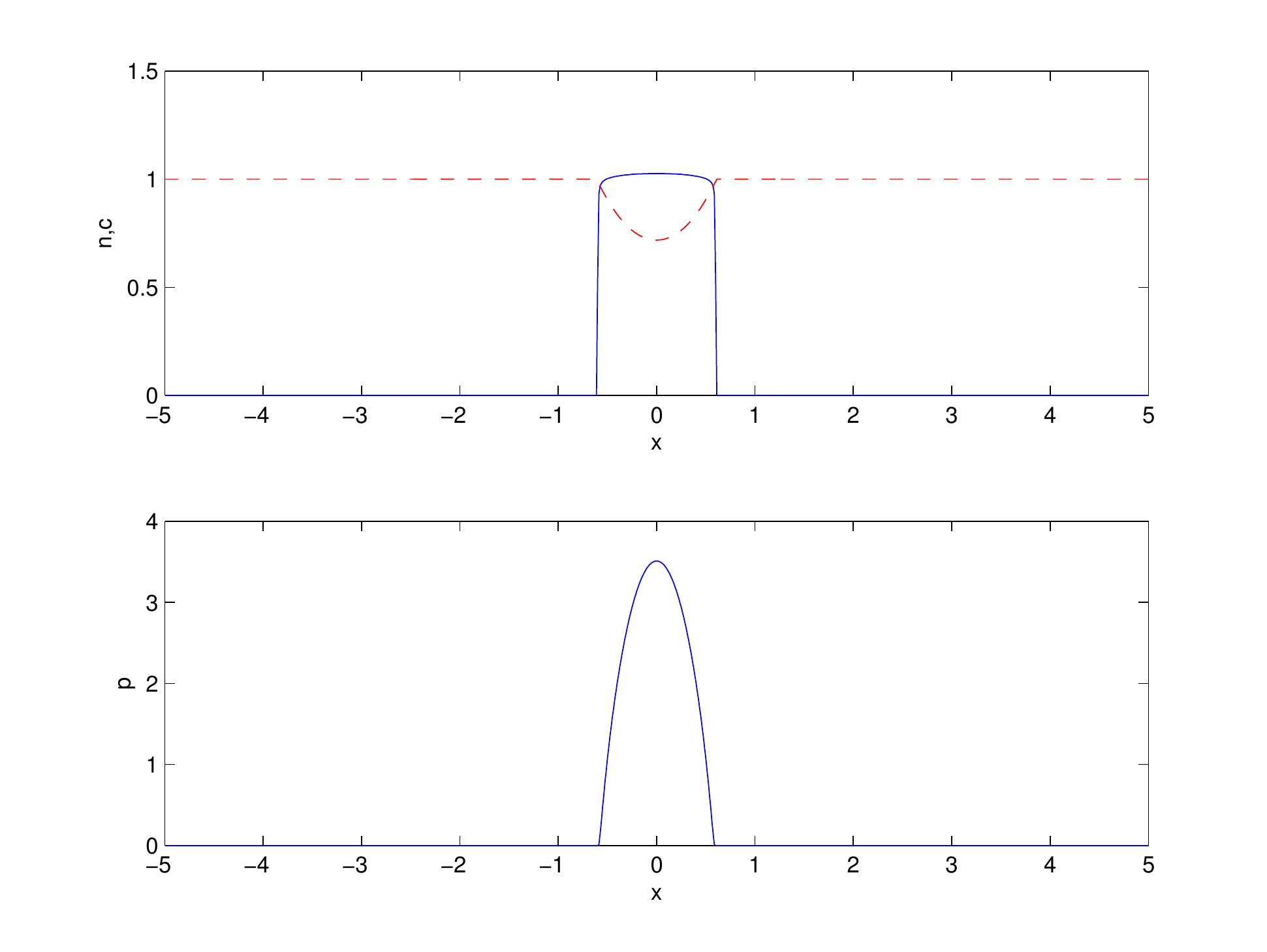}
c)\includegraphics[width=0.45\textwidth]{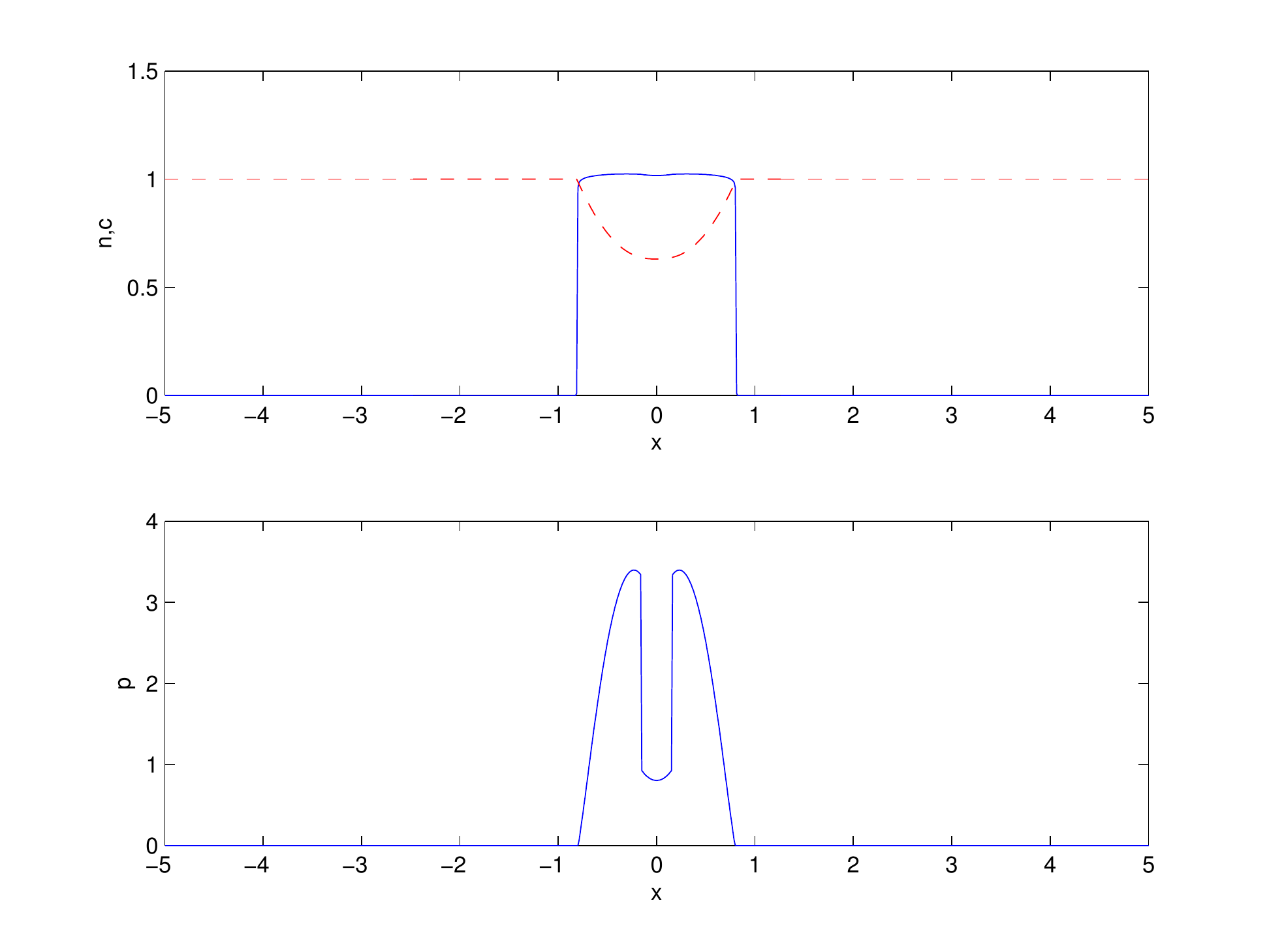}
d)\includegraphics[width=0.45\textwidth]{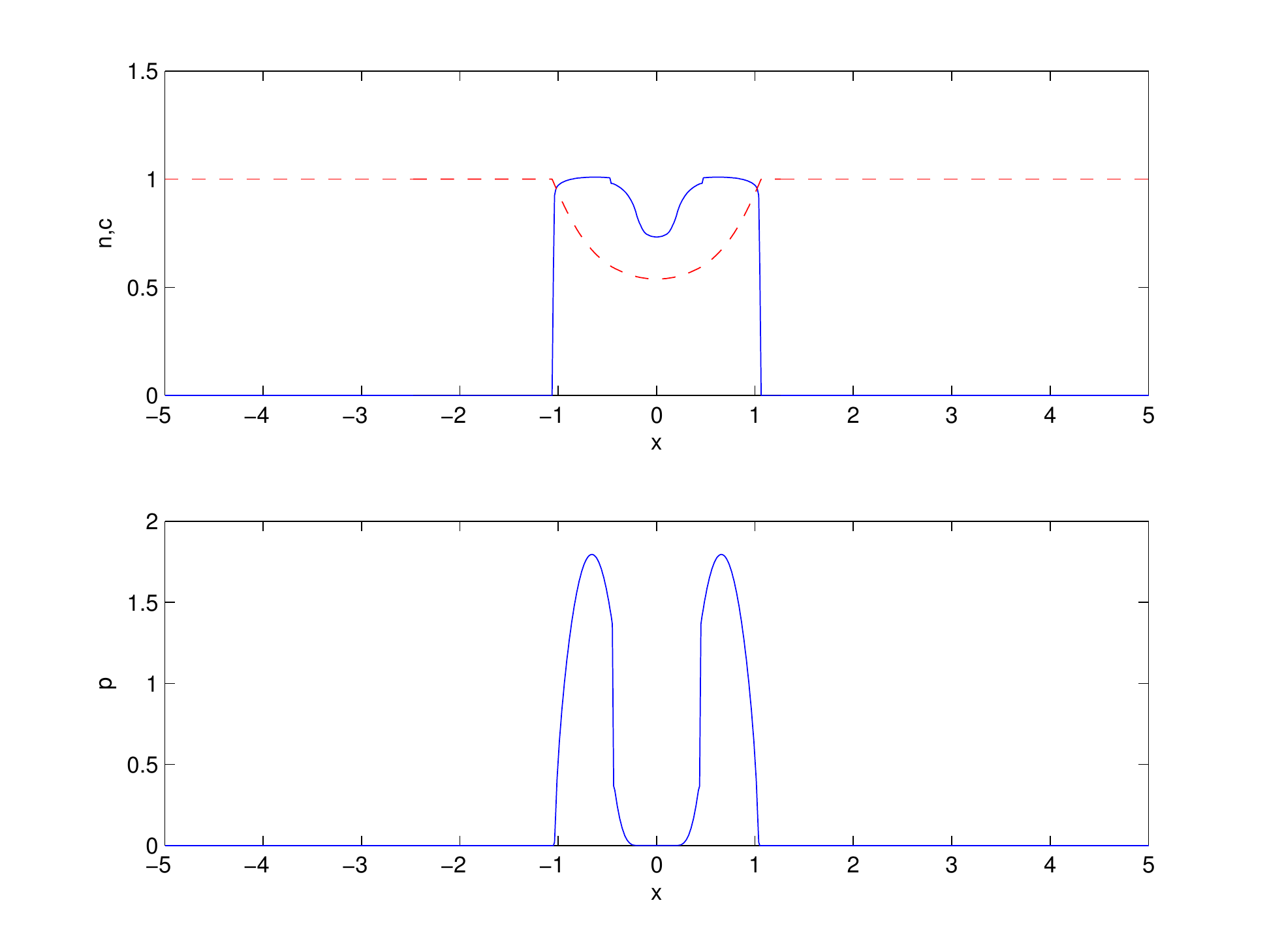}
e)\includegraphics[width=0.45\textwidth]{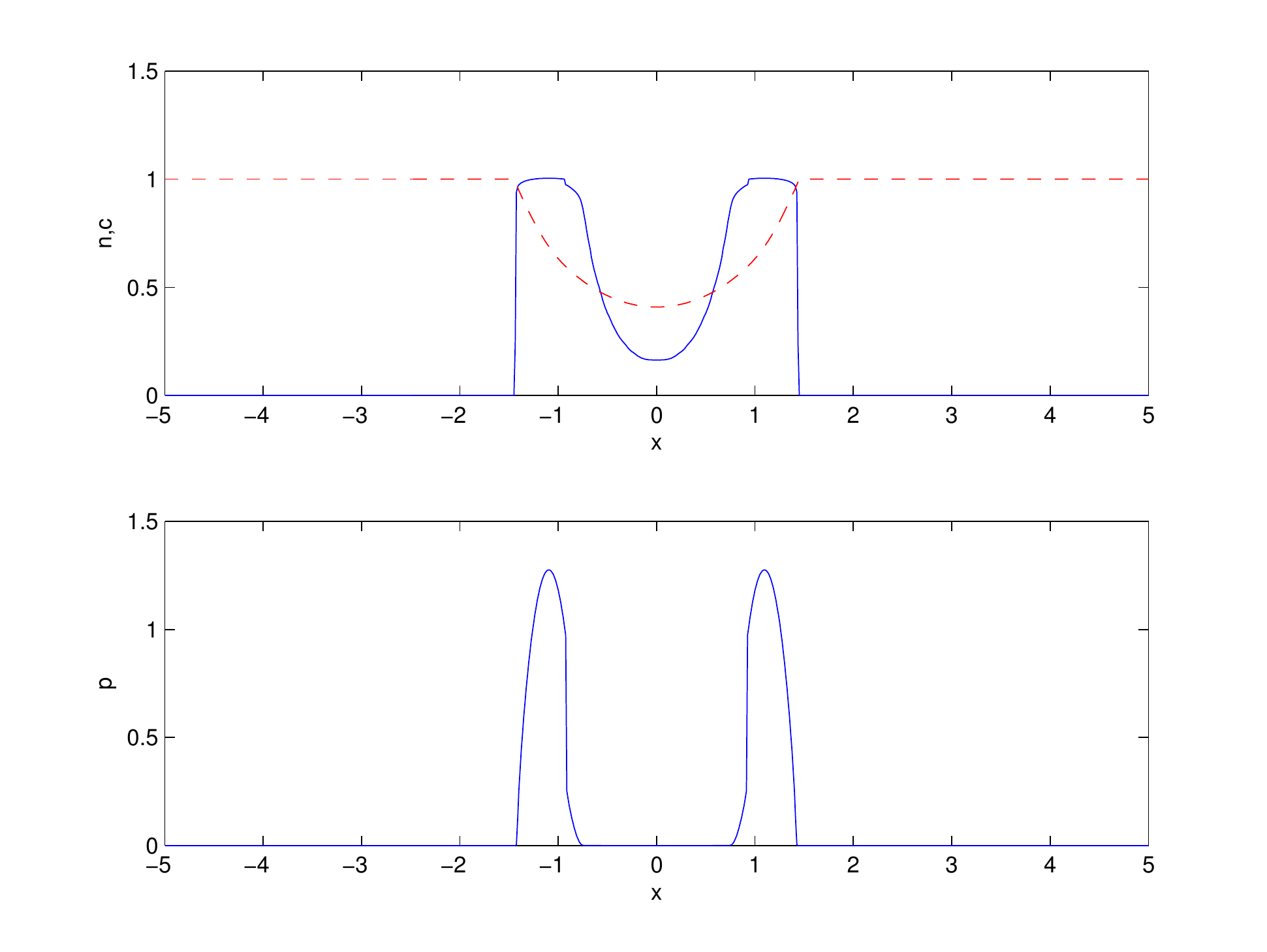}
f)\includegraphics[width=0.45\textwidth]{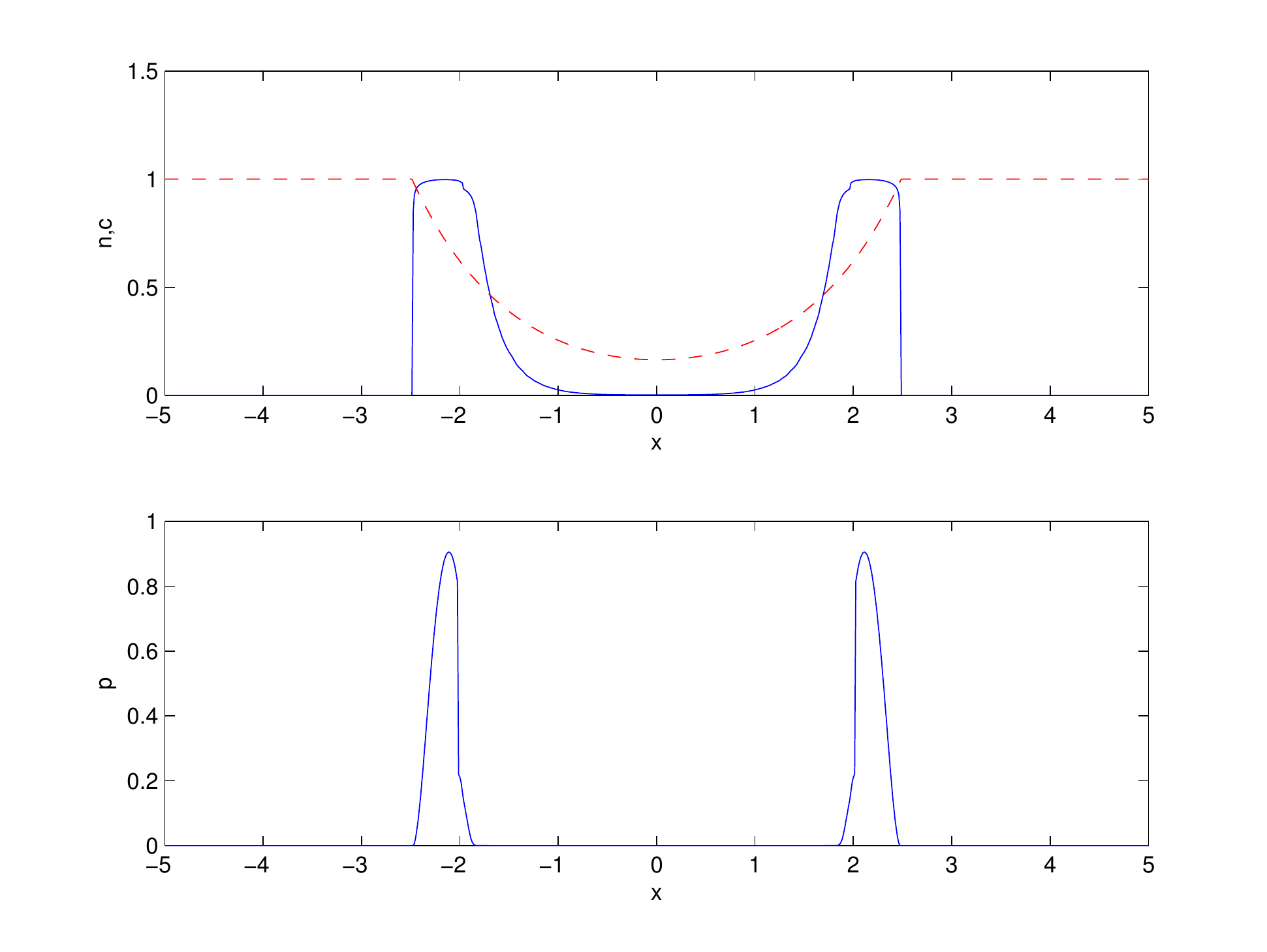}
\caption{Time dynamics of the system \eqref{eq:pqn}.
The parameters are chosen as in \eqref{eq:ex1parameter}.
a), b), c), d), e), f) are for different time and the time increases from a) to f). In each figure, the top subplot depicts
the density distribution $n(x)$ (solid line) and nutrient distribution $c(x)$ (dashed line), 
while the bottom subplot gives the pressure $p(x)$.}
\label{fig:timeevol}
\end{center}
\end{figure}

The numerical solutions are displayed in Fig.~\ref{fig:timeevol} for six different times. We can see that initially the
cell density increases without motion, see subfig. a) and b). After a while, the pressure increases and
the front of the tumor begins to move outwards, see subfig. b) and later. As the size of the tumor increases, the cell
density in the middle decreases due to the lack of nutrient, see subfig. d) and later. In this one dimensional computation, 
we finally obtain two plateaus  that separate, with tails that enlarge and apparently let a necrotic core appear in the very center (the analytic form show that the cell number density does not really vanish, it is in fact an exponential decay). During all these
process, $c$ is constant outside of the tumor as expected.
\\

\begin{figure}
\begin{center}
a)\includegraphics[width=0.45\textwidth]{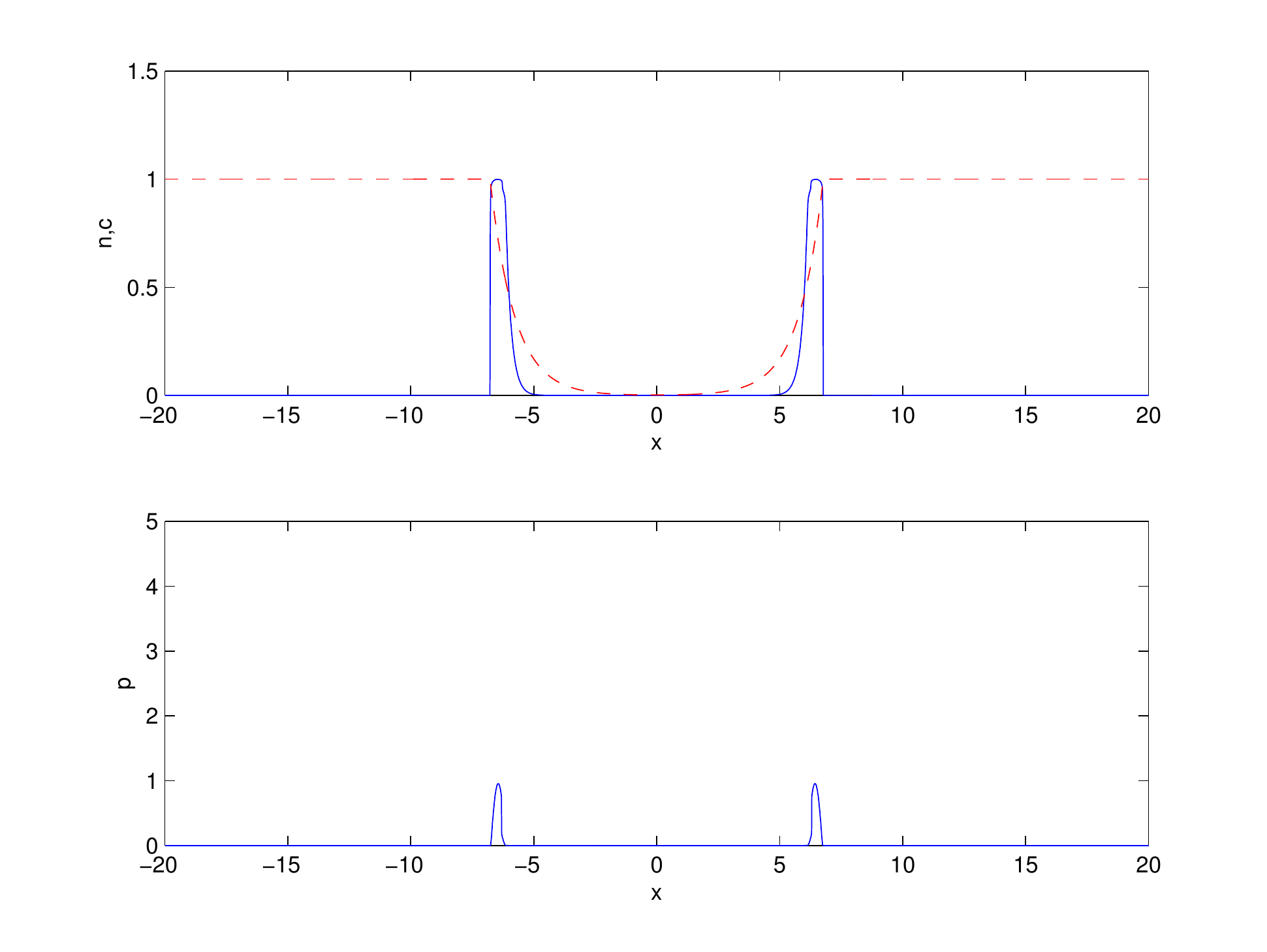}
b)\includegraphics[width=0.45\textwidth]{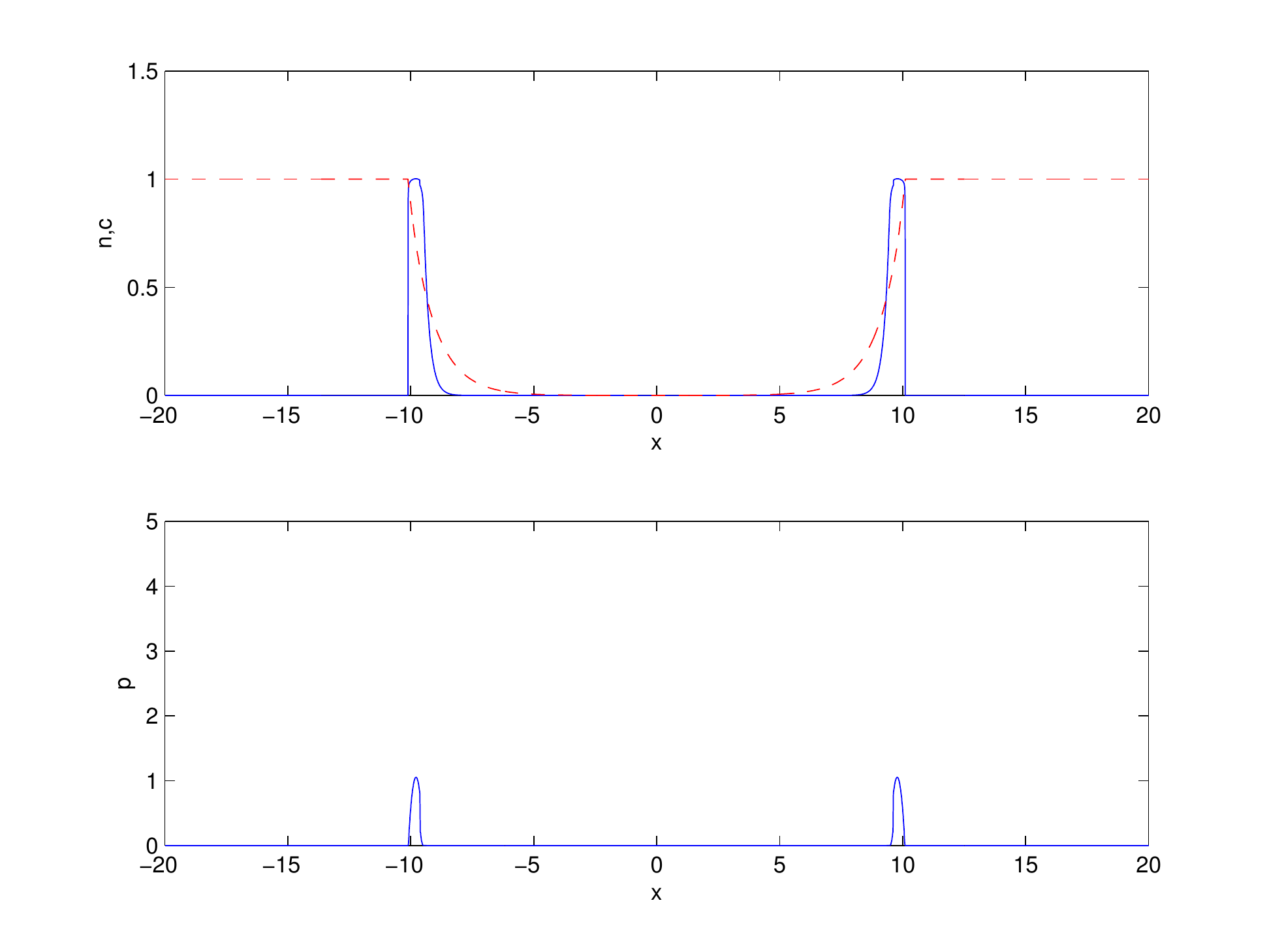}
c)\includegraphics[width=0.45\textwidth]{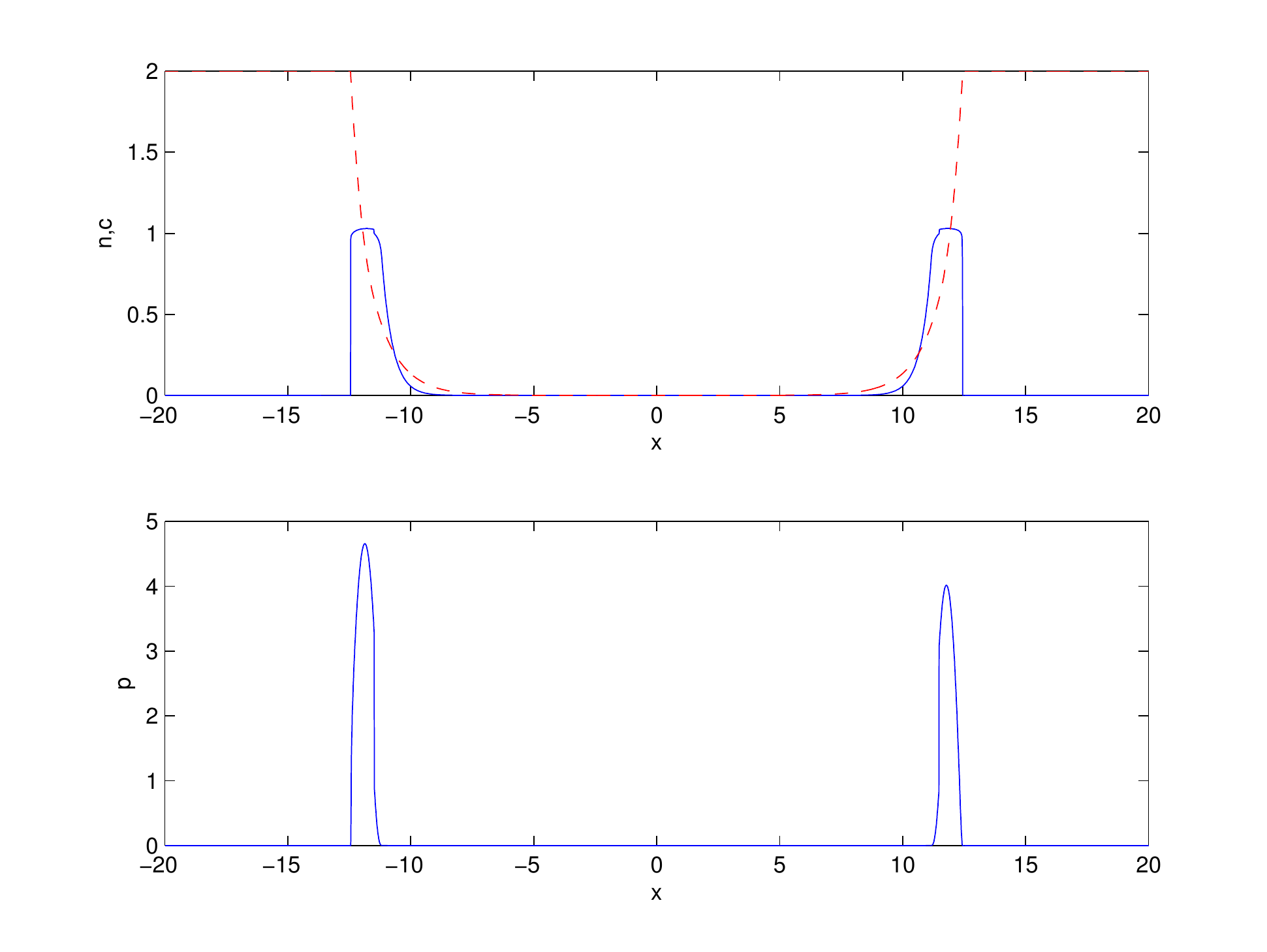}
d)\includegraphics[width=0.45\textwidth]{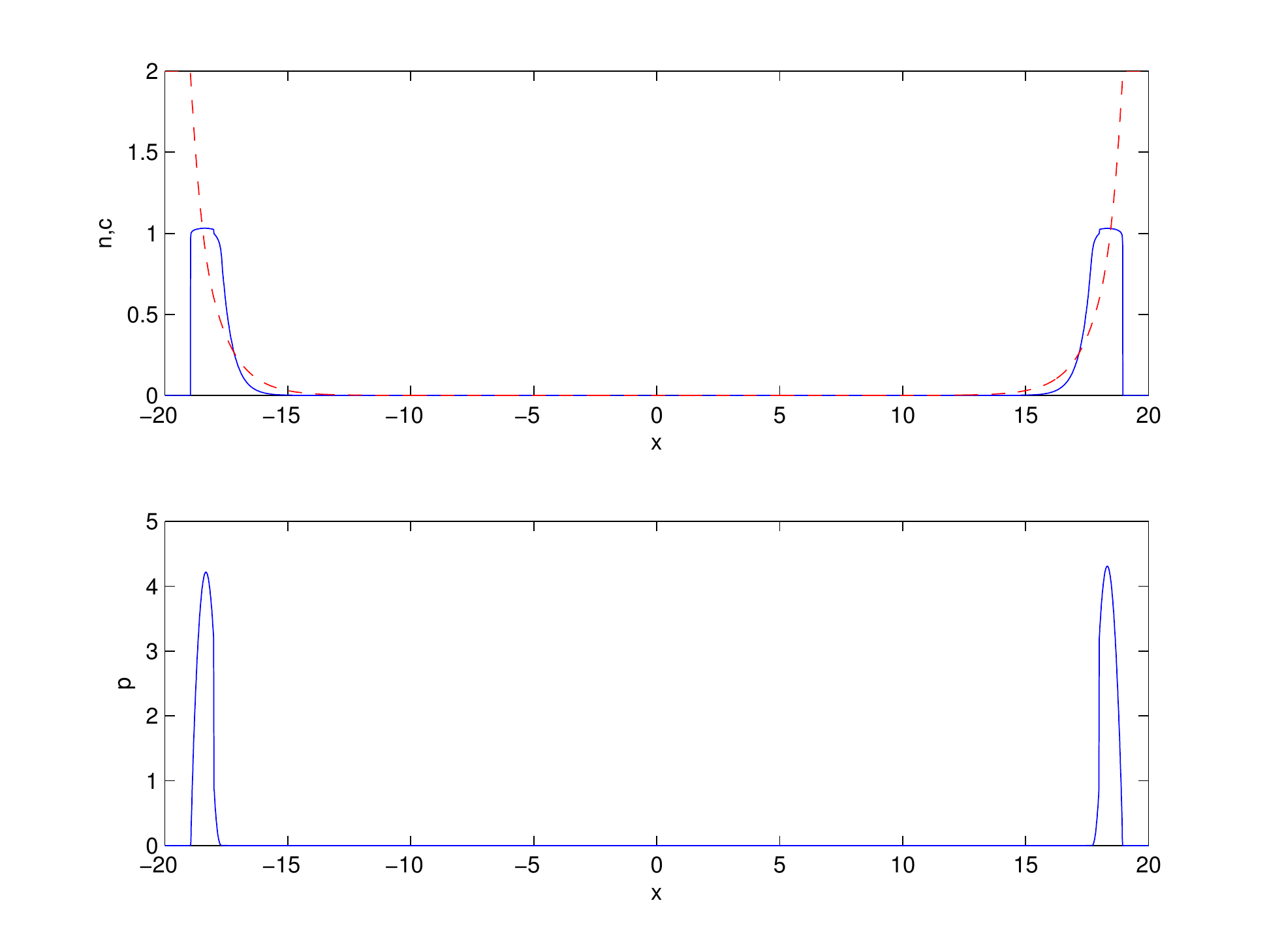}
\caption{Time dynamics {\it in vitro} computed with a larger domain $[-20,20]$.
a), b) are for $c_B=1$ and c), d) are for $c_B=2$. Here a), d) are the 
results at time $t=1$, b), e) are at $t=1.5$.  In each figure, the top subplot depicts
the density distribution $n(x)$ (solid line) and nutrient distribution $c(x)$ (dashed line), 
while the bottom subplot gives the pressure $p(x)$.}
\label{fig:ex1_TW}
\end{center}
\end{figure}

In order to see the traveling wave solutions, we can use a larger computational domain $[-20,20]$
with the same parameters and initial condition.
The results are displayed in Fig.~\ref{fig:ex1_TW}. We can observe numerically that the traveling velocities 
and the sizes of the plateaus increase with $c_B$. Besides, the maximum of the pressure for $c_B=2$ is almost
four times  its maximum value for $c_B=1$.
\\
%

\begin{figure}
\begin{center}
a)\includegraphics[width=0.45\textwidth]{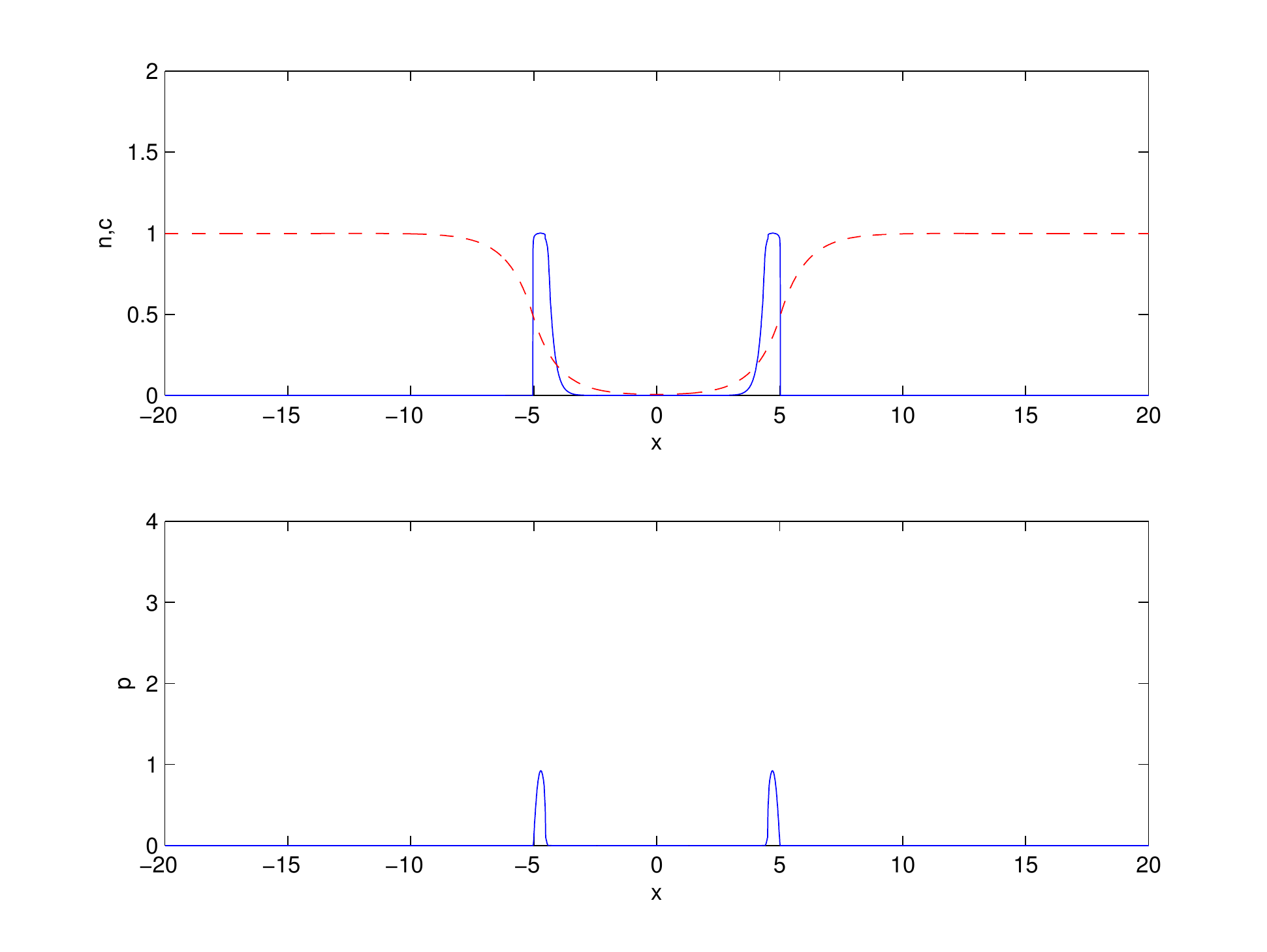}
b)\includegraphics[width=0.45\textwidth]{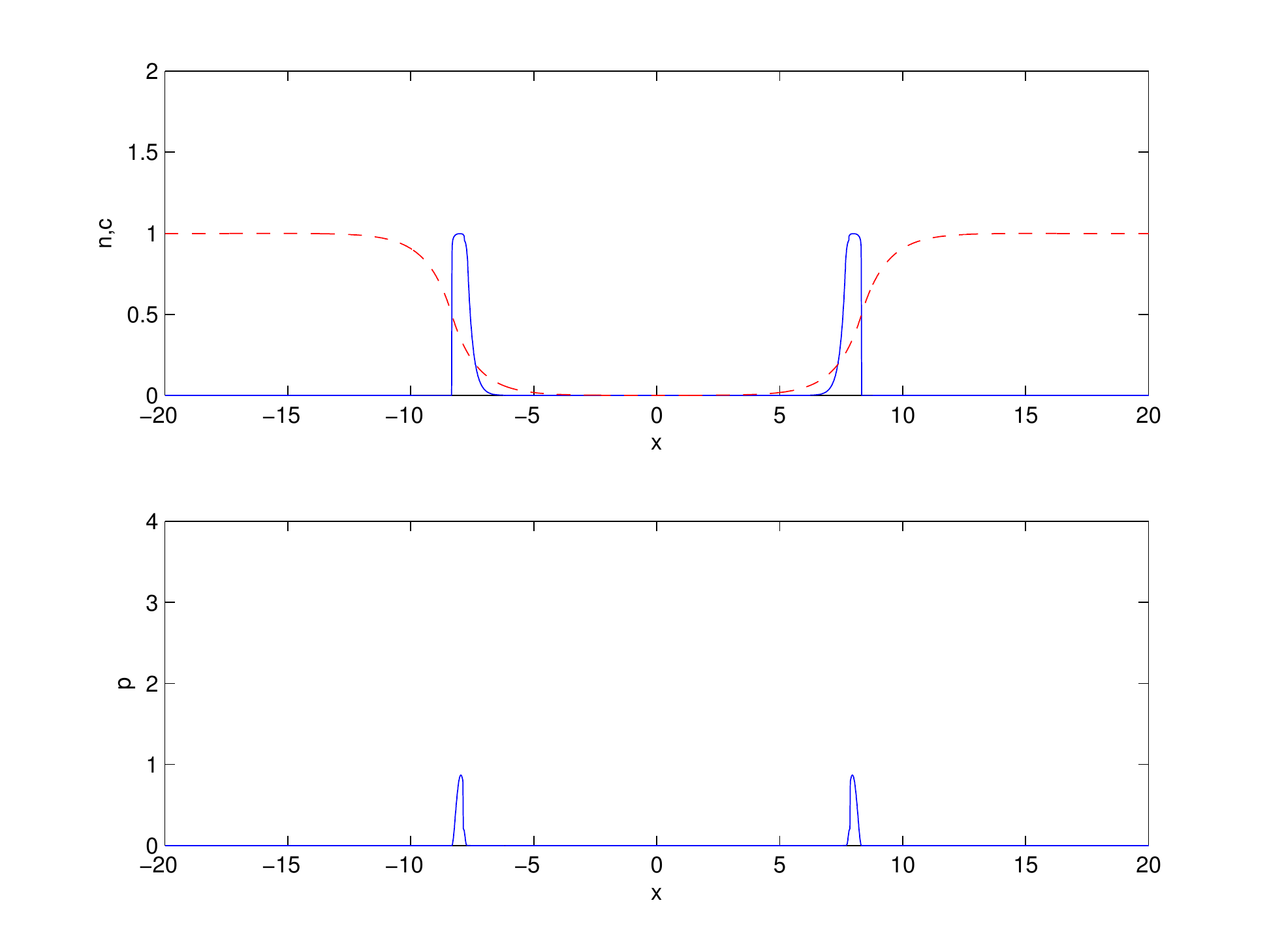}
c)\includegraphics[width=0.45\textwidth]{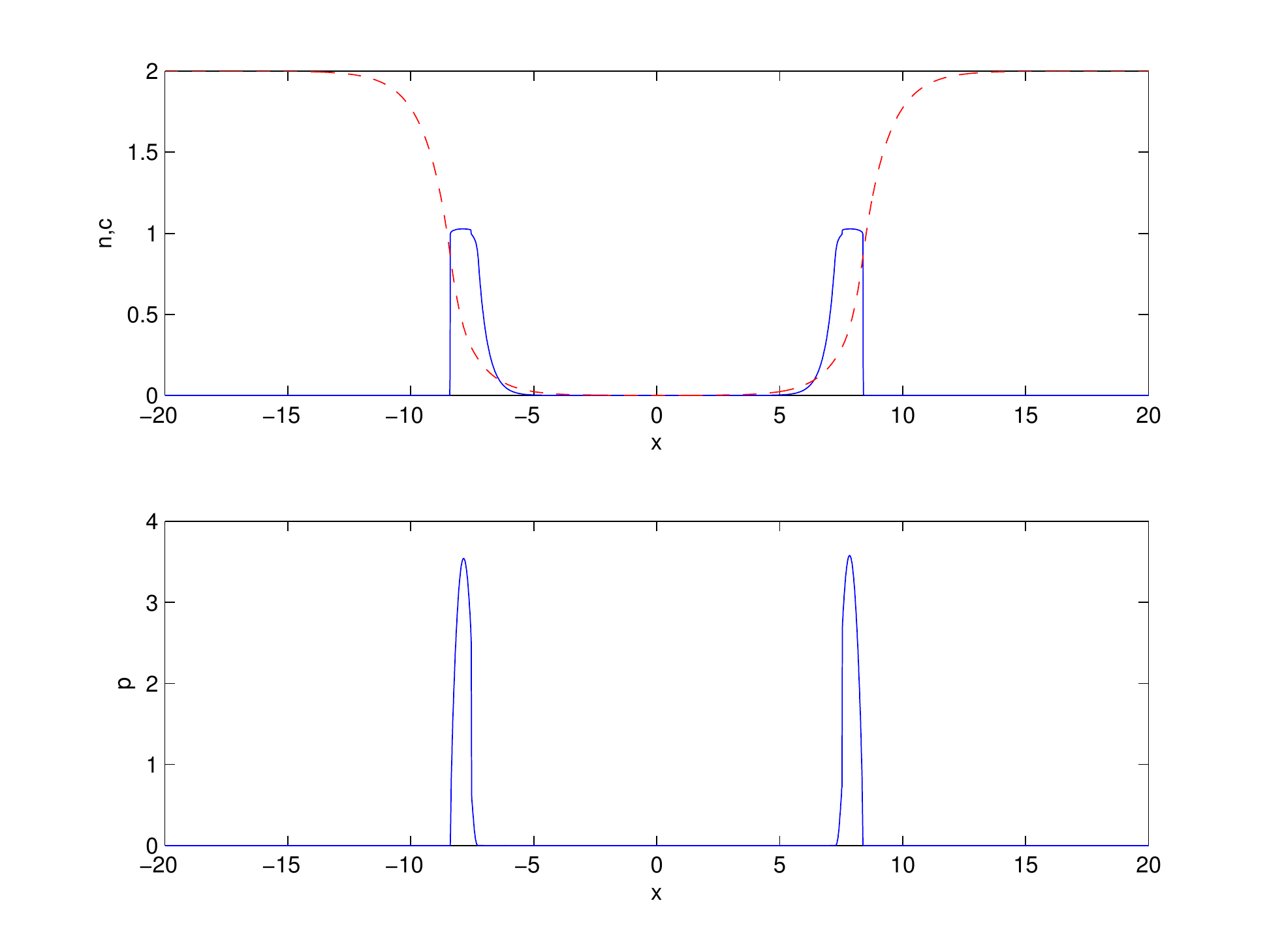}
d)\includegraphics[width=0.45\textwidth]{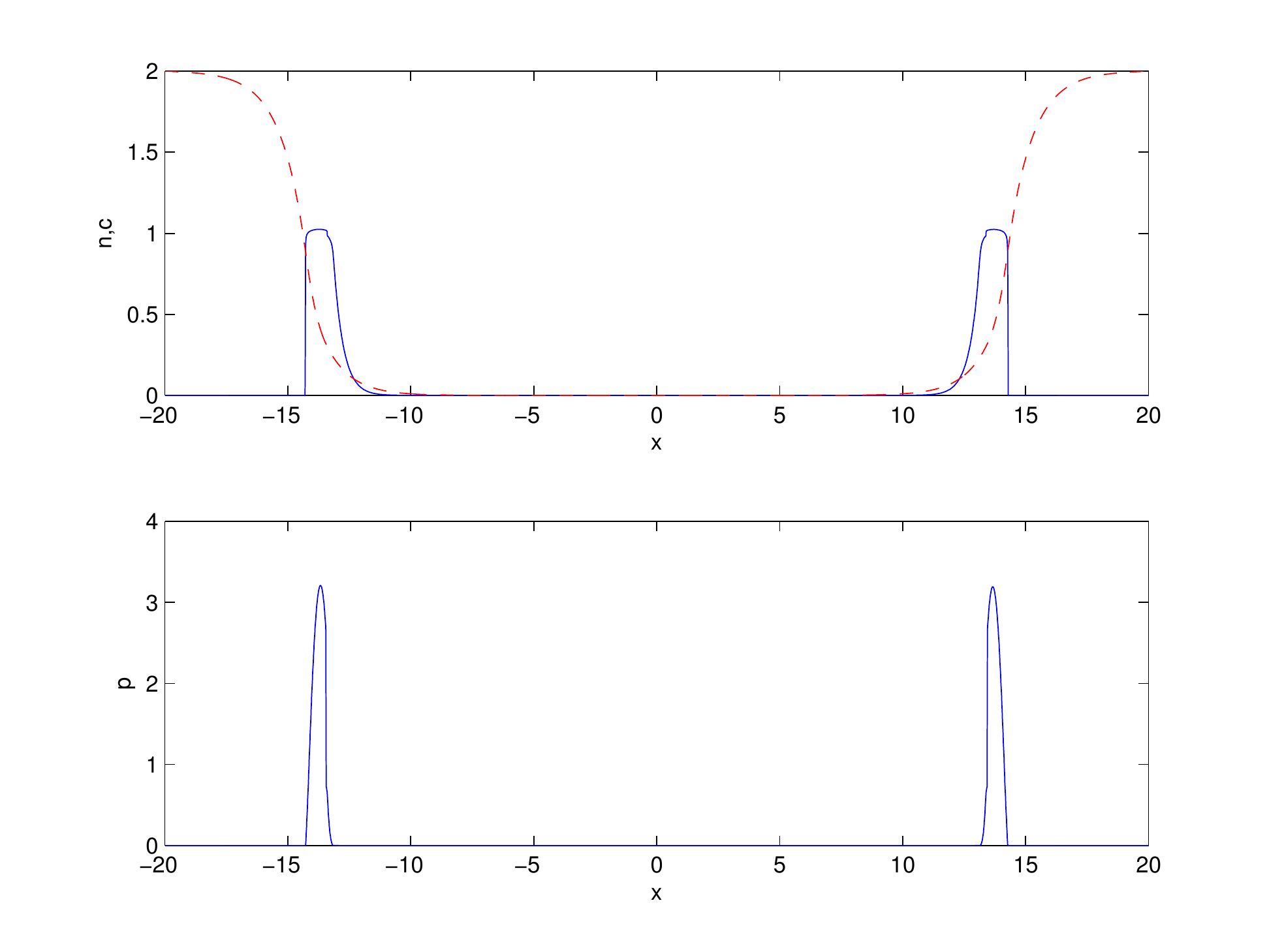}
\caption{Time dynamics {\it in vivo}  computed with a larger domain $[-20,20]$.
a), b) are for $c_B=1$ and c), d) are for $c_B=2$. Here a), c) are the 
results at $T=0.75$ and b), d) are at $T=1$.
In each figure, the top subplot depicts
the cell distribution $n(x)$ (solid line) and nutrient concentration $c(x)$ (dashed line), 
while the bottom subplot gives the pressure $p(x)$.}
\label{fig:ex1_TWv}
\end{center}
\end{figure}

{\it In vivo}, that is for the equation \eqref{eq:vitroc},  the traveling waves are displayed in Fig.~\ref{fig:ex1_TWv}. Here, we use the same values of $\gamma$ and $\psi(n)$ as in \eqref{eq:ex1parameter} but the growth function $G(c)$ is given by 
\beq\label{num:Ginvivo}
G(c)=\left\{\begin{array}{ll}
\ \  21,&c\geq 0.3,\\
-30,&c<0.3.\end{array}\right.
\eeq
With the same initial and boundary conditions as for the {\it in vitro} case,
the nutrient decreases exponentially from outside of the tumor and wave velocity increases with $c_B$.
If we keep using the $G(c)$ as in \eqref{eq:ex1parameter}, the initial cell density decreases to zero as time goes on. 

%

\begin{figure}
\begin{center}
a)\includegraphics[width=0.45\textwidth]{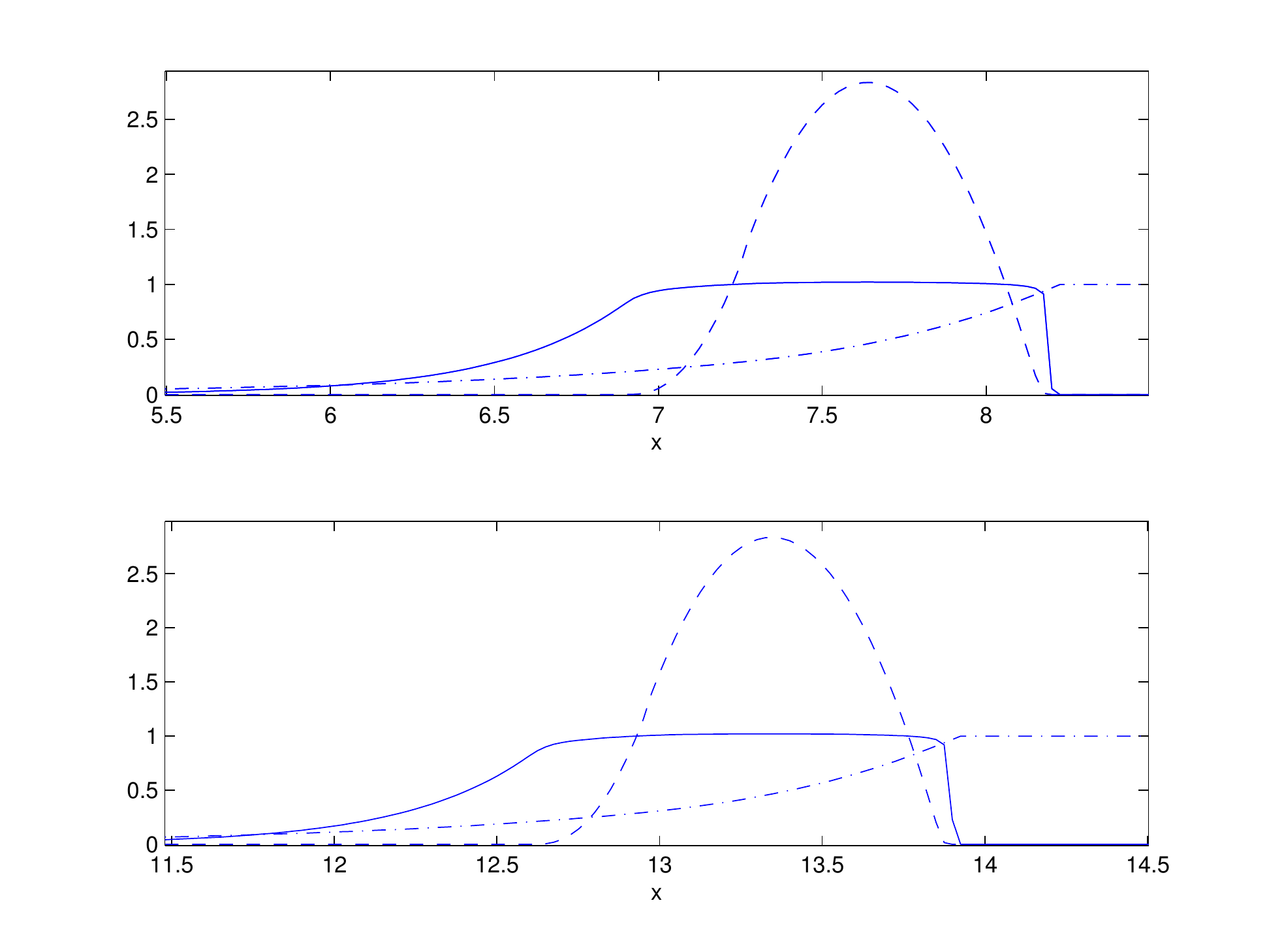}
b)\includegraphics[width=0.45\textwidth]{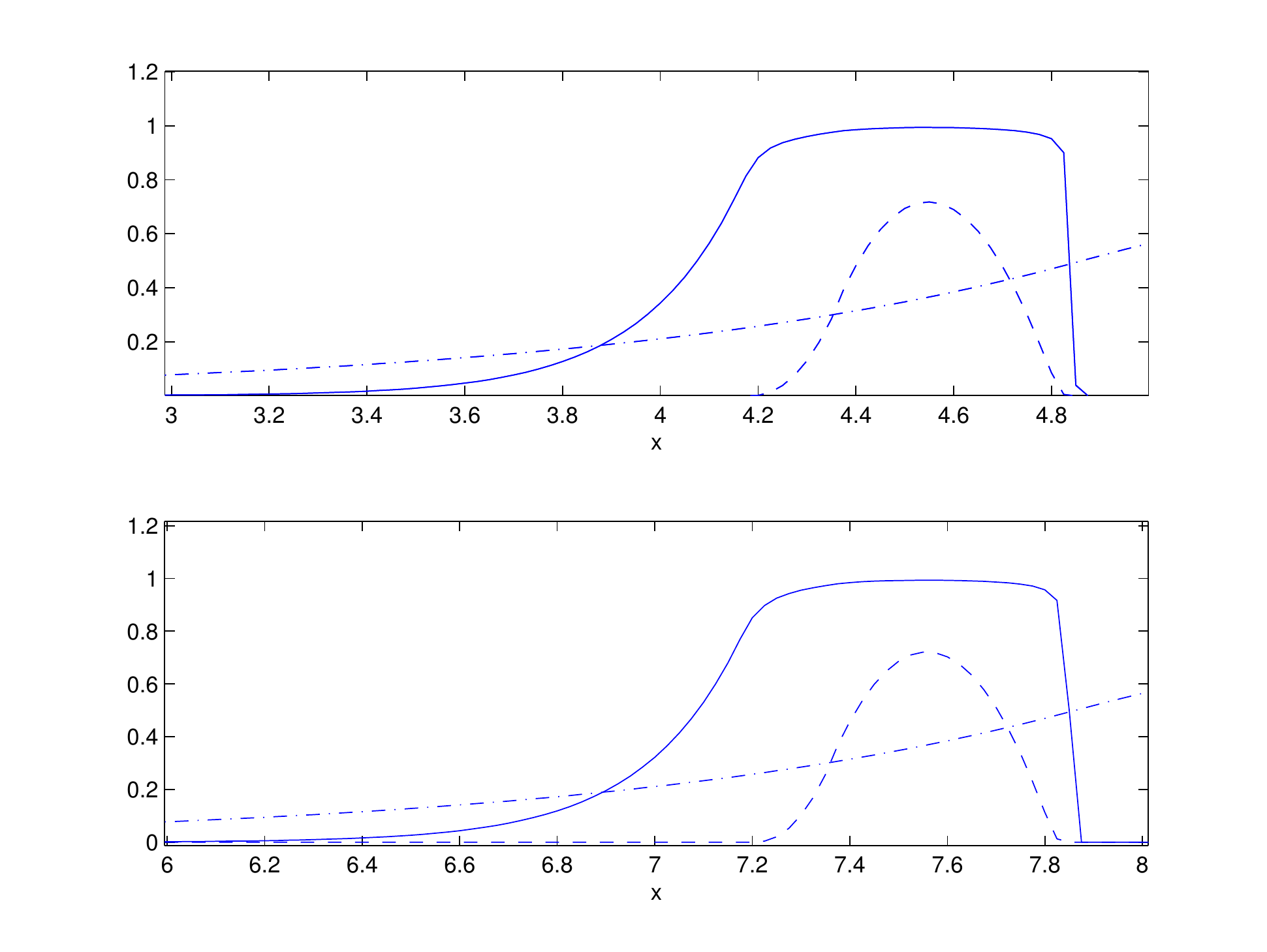}
\caption{Zoom of the pressure profiles for both {\it in vitro} (a) and {\it in vivo} (b) models. Results are displayed at $T=0.75$ (top) and $T=1.25$ (bottom)
with the same parameters and initial conditions.
Each figure depicts
cell distribution $n(x)$ (solid line), nutrient concentration $c(x)$ (dash dotted line), 
and pressure $p(x)$ (dashed line).}
\label{fig:ex1_TWcomp}
\end{center}
\end{figure}

The pressure profiles for both {\it in vitro} and {\it in vivo} models
is displayed in Fig.~\ref{fig:ex1_TWcomp} for $c_B=1$
and with the same parameters and initial conditions.
Smooth transition of the pressure to the necrotic core can be observed in both cases.
Besides the traveling velocity is a bit larger than
the one computed analytically in the subsequent part of this paper.
This is because additional numerical diffusion is introduced at the front
of the pressure.
Such phenomena is similar as in simulations of the porous media equation.

\section{Traveling waves for a simplified Hele-Shaw model with nutrients}
\label{sec:TW}

Traveling waves correspond to an established regime and are a convenient way to explain patterns in cancer invasions \cite{BHMW,TVCVDP}. Therefore, as observed numerically in the previous section, we can expect existence of traveling waves with a proliferative rim and a necrotic core that
invades the domain, at least for large values of $\gamma$ in equation \eqref{eq:pqn} coupled to \eqref{eq:vitroc} or \eqref{eq:vivoc}. It is not our purpose here, to give a  general existence proof of such traveling waves for a fixed and
finite $\gamma$. Here, we focus on the asymptotic model, which is the following
free boundary Hele-Shaw model with nutrients~:
\beq\label{eq:HSM}
\left\{
\begin{array}{ll}
\dis \p_tn-\dv(n\nabla p)=n G(c), \qquad
&\mbox {in $\Omega_0(t)$}, \qquad \Omega_0(t)=\{p(t,x)=0\},
\\[1mm]
n(x,t)=1&\mbox{in $\Omega(t)$},\qquad \ \  \Omega(t) \; =\{p(t,x)>0\},
\\[1mm]
-\Delta c+ \Psi(n,c)=0, &\dis \lim_{x\to + \infty} c(x) = c_B, 
\\[1mm]
-\Delta p=G(c),& \mbox{in $\Omega(t)$},
\\[1mm]
p=0,& \mbox{on $\p \Omega(t)$}.
\end{array}
\right.
\eeq
The growth term $G$ is assumed to satisfy the conditions \eqref{as:tgnbis}.
\\

Traveling waves are solutions which can be written under the form 
$$
n(t,x)=n(x\cdot\mathbf{e}-\sigma t), \quad c(t,x)=c(x\cdot\mathbf{e}-\sigma t), \quad p(t,x)=n(x\cdot\mathbf{e}-\sigma t),
$$
where $\sigma>0$ is a constant representing the traveling wave velocity.
This leads to the following, time independent,  system~:
\beq\label{eq:HSMtwmd}
\left\{
\begin{array}{ll}
-\sigma \mathbf{e}\cdot\nabla n =n G(c),\qquad 
& \mbox{in $\Omega_0(t)$}, \qquad \Omega_0=\{p(x)=0\},
\\[1mm]
n=1,&\mbox{in $\Omega(t)$},\qquad \ \  \Omega \; =\{p(x)>0\},
\\[1mm]
-\Delta c+ \Psi(n,c)=0,& \dis \lim_{x\to +\infty} c(x) = c_B,
\\ [1mm]
-\Delta p=G(c),& \mbox{in $\Omega(t)$},
\\[1mm]
p=0,& \mbox{on $\p \Omega(t)$}.
\end{array}
\right.
\eeq
For the two different models  under consideration, {\it in vitro} and {\it in vivo}, the general equation for $c$ in \eqref{eq:HSMtwmd} is defined by  \eqref{eq:vivoc} and \eqref{eq:vitroc} respectively.

\subsection{Existence of traveling waves in one dimension {\it in vitro}}

Focusing on the one dimensional case  {\it in vitro}, we look for a traveling wave on the
real line propagating to the left. Then, we can simplify the system \eqref{eq:HSMtwmd}
 by introducing the parameter $R>0$ such that 
$$
\Omega=\{p(x)>0\}=[0,R], \qquad \qquad \{n(x)>0\}=(-\infty,R).
$$ 
The above system becomes
\begin{eqnarray}
& -\sigma n' = nG(c), \quad \mbox{ in } (-\infty,0], \label{eq:TWg1}
\\
& n  =1 \quad \mbox{ on } [0,R], \qquad  n=0 \ \mbox{ on } (R,+\infty), \label{eq:TWg2}
\end{eqnarray}
\begin{equation}
\label{eq:TWg3}
 c''  = \psi(n) c \ \mbox{ on } (-\infty,R], \qquad  c(R) = c_B ,
\end{equation}
\begin{equation}
\label{eq:TWg4}
 - p'' = G(c)  \quad \mbox{ on } [0,R], \qquad  p(0)=p(R)=0,  \qquad p \geq 0.
\end{equation}
We claim, this can be completed with the following jump relations at $x=0$ and $R$,
\begin{equation}
  \label{eq:jumprel}
\sigma = -p'(R^-), \qquad n(0) = 1, \qquad p'(0^+)=0.
\end{equation}

Indeed, together with \eqref{eq:TWg4}, the jump relations at $x=0$ and $R$, for the equation on $n$, give
$[\sigma n]+[n\p_x p]=0$. Therefore, denoting $n_0=n(0^-)\in (0,1]$, we arrive at
$$
p'(0^+)=-\sigma(1-n_0),\qquad p'(R^-)=-\sigma, \qquad n_0=n(0^-)\in (0,1].
$$
Since $p$ is constructed as the asymptotic limit of $n^\gamma$ when 
$\gamma\to+\infty$, we restrict our study to nonnegative pressures.
However, at the point $0$, we have $p(0)=0$, therefore the derivative 
should satisfy $p'(0^+)\geq 0$, in order to have $p(x)$ nonnegative when $x>0$
in the neighborhood of $0$. Since $n_0\leq 1$, from the above relation,  we deduce \eqref{eq:jumprel}.

Since we wish to build the traveling waves with semi-explicit formulas, we impose additional conditions for the $C^1$ functions $G$ and $\psi$~:
\begin{equation}  \label{hypG}
  G'\geq 0, \quad \exists\, \bar{c}>0 \mbox{ such that } G(c)=-g_-<0 \mbox{ for } c<\bar{c}, \mbox{ and }
G(c)>0 \mbox{ for } c>\bar{c},
\end{equation}
\begin{equation}
  \label{hyppsi}
\forall\, z\in (0,1), \ \  0<\psi( z )\leq \psi(1), \quad  \psi(0)=0.
\end{equation}
This latter assumption is automatically satisfied if 
$\psi$ is nonincreasing. Also, since $\psi\in C^1([0,1])$ and $\psi(0)=0$, it implies in particular that $\int_0^1\f{\psi(z)}{z}\,dz < \infty$, a property that is used later on.

With these asumptions, we can state an existence result of a traveling wave solution for the 
Hele-Shaw system with nutrients
\begin{theorem}\label{th:TWexist}
Let us assume that $c_B>\bar{c}>0$ and that $G$ and $\psi$
are $C^1$ functions satisfying assumptions \eqref{hypG}-\eqref{hyppsi}.
Then, there exist $\sigma>0$ and $R>0$ such that the system 
\eqref{eq:TWg1}--\eqref{eq:jumprel} admits a solution with $c$ 
increasing on $(-\infty,R]$, 
$n$ increasing on $(-\infty,0]$ and $\dis \lim_{x\to -\infty} n(x) = 0$.
\end{theorem}

\noindent {\bf Proof.} The idea of the proof is as follows. We give a value $\sg>0$ and build a solution $(n_\sg,c_\sg, p_\sg)$ of \eqref{eq:TWg1}--\eqref{eq:jumprel} just discarding the relation for $\sg$ in \eqref{eq:jumprel}. To do so, we need to  fix the parameter $R_\sg$.  Therefore, we split the line into the necrotic, the proliferative and the healthy regions
$$
(0, \infty)= I\cup II \cup III=(-\infty,0)\cup [0,R_\sg]\cup (R_\sg,+\infty),
$$
and we build the solution successively on each interval. 

Then, the value of $\sg$ is determined by the fixed point problem
\beq
\sg = - p_\sg'(R_\sg).
\label{sp:fixedpoint}
\eeq
\noindent {\em Step 1. The piecewise construction of the wave.} Fig. \ref{fig:ex1_TWcomp} can serve as a guide to vizualize the construction which follows. We build the solution departing from the last (healthy) region. Indeed, we have
\\[5pt]
$\bullet$ On $III=(R,+\infty)$, $n(x)=0$, $p(x)=0$, $c=c_B$.
\\[5pt]
$\bullet$
On $II=[0,R]$, $n(x)=1$. Then, the system reduces to the to solving the
equations
$$
-c'' + \psi(1) c = 0, \qquad c(R)=c_B,
$$
$$ 
-p'' = G(c), \qquad p(0)=p'(0)=0.
$$
Setting $c'_R=c'(R)/\sqrt{\psi(1)}$ which is unknown, 
we obtain by solving the first equation
\begin{equation}
  \label{cII}
c(x)=c_B \cosh\big(\sqrt{\psi(1)}(x-R) \big) + c'_R \sinh\big(\sqrt{\psi(1)}(x-R)\big).
\end{equation}
For the second equation, we have
\begin{equation}
  \label{pII}
p(x) = -\int_0^x\int_0^y G(c)(z)\,dzdy = -\int_0^x (x-z) G(c)(z)\,dz.
\end{equation}
Moreover, the boundary condition $p(R)=0$, gives
$$
\int_0^R\int_0^y G\big(c_B \cosh(\sqrt{\psi(1)}(z-R)) + 
c'_R \sinh(\sqrt{\psi(1)}(z-R))\big) \,dzdy = 0.
$$
Applying the Fubini theorem, we can rewrite this last equation as
$$
\int_0^R z\,G\left(c_B \cosh(\sqrt{\psi(1)}z) - 
c'_R \sinh(\sqrt{\psi(1)}z)\right) \,dz = 0.
$$
Setting $s=z/R$, this nonlinear equation gives a first relation between the two free parameters $R$ and $c'_R$
\begin{equation}
  \label{eq1}
\int_0^1 s \; G\left( c_B\cosh \big(\sqrt{\psi(1)}Rs \big)- c'_R \sinh\big(\sqrt{\psi(1)}Rs \big)\right) \,ds = 0.
\end{equation}

Moreover, we are looking for a solution $c$ which is positive
and nondecreasing. From \eqref{cII}, we have
$$
c'(x)=c_B \sqrt{\psi(1)} \sinh\big(\sqrt{\psi(1)}(x-R) \big) + 
c'_R \sqrt{\psi(1)} \cosh\big(\sqrt{\psi(1)}(x-R)\big).
$$
Thus, $c$ is positive and nondecreasing on $[0,R]$ iff $c(0)>0$ and $c'(0)>0$, that is 
\begin{equation}
  \label{boundcR}
c_B \tanh(\sqrt{\psi(1)}R) \leq c'_R < \frac{c_B}{\tanh(\sqrt{\psi(1)}R)}.
\end{equation}
Notice that, with \eqref{boundcR} and \eqref{eq1}, we have necessarily
\begin{equation}
  \label{hypc0}
\bar{c}>c(0 ) > 0 \quad \mbox{ and } \quad c' (0) \geq 0.
\end{equation}
because $G(c)$ has to be negative on a subinterval of $(0,R)$,
and  monotonicity of $G$ and $c$.
\\
\\
$\bullet$  On $I=(-\infty,0)$, we have $p=0$ and we look for a 
solution to the system
\begin{eqnarray}
  \label{eqg:nI}
&  -\sigma n' = n G(c), \qquad n(0)=1, \\
  \label{eqg:cI}
&  c'' = \psi(n) c, \qquad c(0)=c_0, \ c'(0)=c'_0 ,
\end{eqnarray}
where by continuity for $c$ and $c'$ at $0$,  from \eqref{cII} we conclude that
\begin{eqnarray}
&\dis  \label{eq:c0}
c_0=c_B \cosh(\sqrt{\psi(1)}R) - c'_R \sinh(\sqrt{\psi(1)}R) < \bar c, \\[2mm]
&\dis \label{eq:cprime0}
c'_0=c'_R\sqrt{\psi(1)} \cosh(\sqrt{\psi(1)}R) - c_B \sqrt{\psi(1)}\sinh(\sqrt{\psi(1)}R) >0.
\end{eqnarray}

From \eqref{hypG}, we deduce that $G(c)=-g_-$ on $(-\infty,0]$.
Therefore the solution of equation \eqref{eqg:nI} is given by 
\begin{equation}
  \label{eq:nIexpr}
n(x) = e^{\f{g_-}{\sigma} x}, \qquad \mbox{for } x<0.
\end{equation}
From this, we show in the proof of  Lemma \ref{lem1} below, that we can solve the equation on $c$ and deduce the existence of a function $A(\sg)> 0$, for $\sg>0$,
such that $c'_0=A(\sigma)c_0$ and the solution $c$ of \eqref{eqg:cI} is
nonnegative and nondecreasing for all $x\in I$.
\\

\noindent {\em Step 2. The nonlinear equations for $R_\sg$ and $c'_R$.}
Using this function  $A(\sg)> 0$ and \eqref{eq:c0}-\eqref{eq:cprime0}, we can reduce our construction to the second relation  between $R$ and $c'_R$ (together with \eqref{eq1})
\begin{equation}
  \label{eq:c'R}
c'_R = c_B \frac{A(\sg)\cosh(\sqrt{\psi(1)}R)+\sqrt{\psi(1)}\sinh(\sqrt{\psi(1)}R)}
{A(\sg)\sinh(\sqrt{\psi(1)}R)+\sqrt{\psi(1)}\cosh(\sqrt{\psi(1)}R)}.
\end{equation}
Notice that,  from this expression, a simple but tedious calculation shows 
that \eqref{boundcR} is satisfied.
\\

Finally, for a given $\sg>0$, we have constructed a possible solution on $(-\infty,R_\sg]$ such that
$c$ is nondecreasing and $\lim_{x\to -\infty} c(x)=c_-\geq 0$.
To conclude the construction, we need to prove that there exists $R_\sigma>0$, $c'_R>0$ such that \eqref{eq1} and  \eqref{eq:c'R} are satisfied.

We define, using \eqref{eq:c'R},  
\begin{equation}
  \label{calG}
\begin{array}{ll}
\dis \gamma_\sg(R,s)& :=  \dis c_B\cosh(\sqrt{\psi(1)}Rs) - c'_R\sinh(\sqrt{\psi(1)}Rs) \\[3mm]
&\dis = c_B \frac{A(\sg)\sinh(\sqrt{\psi(1)}R(1-s))
+\sqrt{\psi(1)}\cosh(\sqrt{\psi(1)}R(1-s))}
{A(\sg)\sinh(\sqrt{\psi(1)}R)+\sqrt{\psi(1)}\cosh(\sqrt{\psi(1)}R)}.
\end{array}
\end{equation}
Then, the equation  \eqref{eq1} is also written,
$$
\int_0^1 sG(\gamma_\sg(R_\sg,s)) ds =0.
$$
We are going to show that there is a unique solution $R_\sg >0$, continuous with respect to $\sg$. 
Firstly, from a straightforward computation, \eqref{calG} shows 
that for all $s\in (0,1)$
$$
\begin{array}{ll}
\dis \f{d}{dR}\gamma_\sg(R,s) <
&\dis \f{c_B\sqrt{\psi(1)}\sinh(\sqrt{\psi(1)}R s)\big(A(\sg)^2-\psi(1)\big)}
{\Big(A(\sg)\sinh(\sqrt{\psi(1)}R)+\sqrt{\psi(1)}\cosh(\sqrt{\psi(1)}R)\Big)^2}.
\end{array}
$$
With the estimate of Lemma \ref{lem1}, we deduce that $R\mapsto \gamma_\sg(R,s)$
is a decreasing function for all $s\in (0,1)$.
By monotonicity of $G$ in \eqref{hypG} we deduce that 
$R\mapsto \int_0^1sG(\gamma_\sg(R,s))\,ds$ is a decreasing function.
When $R=0$, we have $\gamma_\sg(0,s)=c_B$, then $\int_0^1 s G(\gamma_\sg(0,s))\,ds = G(c_B)/2 >0$.
When $R\to + \infty$, we have $\gamma_\sg(R,s)\sim e^{-\sqrt{\psi(1)}Rs} \to 0$ for $s>0$;
then $\int_0^1 s G(\gamma_\sg(R,s))\,ds \to  G(0)/2 <0$.
By continuity and monotonicity, we conclude on the existence of 
a unique $R_\sigma>0$ such that \eqref{eq1} is satisfied.
Secondly, by continuity of $\sigma\mapsto A(\sg)$,
we deduce that $\sg\mapsto R_\sg$ is continuous.
\\
\\
%
\noindent {\em Step 3. The nonlinear equation for $\sg$.}
We are reduced to proving the existence of a solution for the fixed point problem
\eqref{sp:fixedpoint}. Thanks to the expression of the pressure, \eqref{pII}, we can rewrite the equation as
$$
\sigma = R_\sigma \int_0^1 G\left(c_B \cosh(\sqrt{\psi(1)}R_\sigma s) - c'_R \sinh(\sqrt{\psi(1)}R_\sigma s)\right) \,ds,
$$
where $c'_R$ is defined by \eqref{eq:c'R}. With the expression of $\gamma_\sg$ given in \eqref{calG}, the latter equation rewrites
\begin{equation}
  \label{eq:sgpt}
  \sigma = R_\sigma \int_0^1 G\big( \gamma_\sg(R_\sigma,s)\big) \,ds.
\end{equation}

We first notice that from Lemma \ref{lem1}, we have $A(\sg)\to 0$ when
$\sg \to 0$. We deduce that $R_\sigma \to R_0>0$ which is a solution to
$$
\int_0^1 sG(\gamma_\sg(R_0,s))\,ds = 
\int_0^1 sG\Big(c_B\f{\cosh(\sqrt{\psi(1)}R_0 (1-s))}{\cosh(\sqrt{\psi(1)}R_0)}
\Big)\,ds = 0.
$$
Since $s\mapsto \gamma_\sg(R_0,s)$ is nonincreasing, there exists
$s_0\in (0,1)$ such that $G(\gamma_\sg(R_0,s))\geq 0$ for all $s<s_0$ and 
$G(\gamma_\sg(R_0,s))<0$ for all $s>s_0$.
We deduce that
$$
0=\int_0^1 sG(\gamma_\sg(R_0,s))\,ds < s_0 \int_0^1 G(\gamma_\sg(R_0,s))\,ds.
$$
Then the limit when $\sg \to 0$ of
the right hand side of \eqref{eq:sgpt} is positive.
On the other hand, from the lower bound of $c'_R$ in \eqref{boundcR}, we have that
$\dis \gamma_\sg(R,s) \leq \frac{c_B}{\cosh(\sqrt{\psi(1)}Rs)}$.
It is straightforward to verify that there is a unique $R_b>0$ such that 
$$
\int_0^1 s G\Big(\frac{c_B}{\cosh(\sqrt{\psi(1)}R_b s)}\Big) \,ds = 0.
$$
We deduce that for all $R\geq R_b$ we have, by monotonicity of $G$,
$$
\int_0^1 s G(\gamma_\sg(R,s))\,ds \leq \int_0^1 sG\Big(\frac{c_B}{\cosh(\sqrt{\psi(1)}Rs)}\Big) \,ds
\leq \int_0^1 s G\Big(\frac{c_B}{\cosh(\sqrt{\psi(1)}R_bs)}\Big) \,ds = 0.
$$
Since $\int_0^1sG(\gamma_\sg(R_\sigma,s))\,ds=0$, we conclude that $0<R_\sigma\leq R_b$.
Thus the 
right hand side of \eqref{eq:sgpt} is bounded.
By continuity, we conclude that equation \eqref{eq:sgpt} admits 
a solution $\sigma>0$.
\qed

\subsection{A technical Lemma}

In the proof of Theorem~\ref{th:TWexist}, we have used an argument which is stated in the 
\begin{lemma}\label{lem1}
Let $\sigma>0$ and $g_->0$ be given and assume that $\psi$
satisfies \eqref{hyppsi}.
We consider the Cauchy problem 
\beq\label{eqlem1}
-c'' + \psi\big(e^{g_-x/\sigma}\big) c = 0 \;  \text{for }  x<0, \qquad c(0)=c_0,\quad c'(0)=c'_ .
\eeq
Then, there exists a continuous mapping $\sg\mapsto A(\sigma)\in (0,\sqrt{\psi(1)} \, ]$ 
such that the solution of the Cauchy problem \eqref{eqlem1} is defined, nonnegative and nondecreasing, if and only if  $c'_0=A(\sigma)c_0$. 
Moreover, the estimate holds~: $A(\sg)\leq  \frac{\sigma}{g_-}\int_{0}^{1} \frac{\psi(\nu)}{\nu}\,d\nu$ and $A(\cdot)$ can be  extended by continuity to $A(0)=0$.
\end{lemma}
{\bf Proof.}
We first reduce equation \eqref{eqlem1} to the following
ODE system,
\begin{eqnarray}
\label{eqc1}& c' = u, \qquad & c(0)=c_0,
\\
\label{equ1}& u' = \psi\big(e^{g_-x/\sigma}\big) c,\qquad & u(0)=c'_0.
\end{eqnarray}
Since we are looking for a nondecreasing $c$, we have $u\geq 0$.
This system can be further simplified by setting 
$y=-e^{g_-x/\sigma}$. We define $\tildeu(y)=u(x)/c_0$ and $\tildec(y)=c(x)/c_0$.
Then system \eqref{eqc1}--\eqref{equ1}, reduces to a Cauchy problem on
$[-1,0]$~:
\begin{equation}
  \label{eq:Cauchypb}
\left\{\begin{array}{ll}
\dis \frac{d\tildec}{dy} = \frac{\sigma}{g_- y}\, \tildeu,
\qquad \frac{d\tildeu}{dy} = \frac{\sigma \psi(-y)}{y \, g_-} \,\tildec, 
\\[3mm]
\dis \tildec(-1)=1, \qquad \quad  \tildeu(-1)=c'_0/c_0.
\end{array}\right.
\end{equation}
There exists a unique solution to this Cauchy problem and $\tildeu$ is nonincreasing and
$\tildec$ is nonincreasing until $\tildeu$ vanishes.
Thus there are two kinds of solutions~:
\begin{itemize}
\item either there exists $y_0\in (-1,0)$ such that $\tildeu(y_0)=0$, this 
solution is called Type I,
\item or $\tildeu>0$ on $(-1,0)$, this solution is called Type II.
\end{itemize}
Obviously, we are looking for a solution in the limiting case, where 
$y_0=0$ and $\tildeu(y)> 0$ for all $y\in[-1,0)$.

We first notice that if $c'_0=0$, then clearly the solution is of Type I.
On the contrary, if we assume that all the solutions are of Type I,
then on $(-1,y_0)$, we have $0\leq \tildec \leq 1$.
It implies from \eqref{eq:Cauchypb} that for $y\in (-1,y_0)$, ($y_0<0$)
$$
\frac{d\tildeu}{dy} = \frac{\sigma \psi(-y) \tildec}{g_-y}
\geq \frac{\sigma \psi(-y)}{g_-y}.
$$
Integrating on $[-1,y_0]$, we deduce that
$$
\tildeu(y_0)\geq \frac{c'_0}{c_0} +\frac{\sigma}{g_-}\int_{-1}^{y_0} \frac{\psi(-\nu)}{\nu}\,d\nu.
$$
The integral in the right hand side is finite thanks to assumption \eqref{hyppsi}.
We deduce that if $c'_0$ is large enough such that
\beq\label{eq:lowerboundAsigma}
\frac{c'_0}{c_0} +\frac{\sigma}{g_-}\int_{-1}^{y_0} \frac{\psi(-\nu)}{\nu}\,d\nu >0,
\eeq
then $\tildeu(y_0)>0$ which contradicts the fact that $\tilde{u}(y_0)=0$.

In summary, for $c'_0=0$ the solution is of Type I, for $c'_0$ large enough
the solution is of Type II. By continuity wich respect to the initial
data, there exists $c'_0$ such that the Cauchy problem \eqref{eq:Cauchypb}
admits a solution such that $y_0=0$ and $\tildeu(y)> 0$ for all $y\in[-1,0)$.

From now on, we denote $A(\sigma)=\tildeu(-1)= \f{c'_0}{c_0}$ to emphasize the
dependency on $\sigma$. Besides, from \eqref{eq:lowerboundAsigma}, 
when $c_0'/c_0>-\frac{\sigma}{g_-}\int_{-1}^{0} \frac{\psi(-\nu)}{\nu}\,d\nu$, the solution is of type II. Thus
we have obviously 
$$
A(\sigma)<-\frac{\sigma}{g_-}\int_{-1}^{0} \frac{\psi(-\nu)}{\nu}\,d\nu=
\frac{\sigma}{g_-}\int_{0}^{1} \frac{\psi(\nu)}{\nu}\,d\nuÊ\quad \forall \sg>0.
$$
Then from the continuity of the Cauchy problem \eqref{eq:Cauchypb} 
with respect to the parameter $\sg$ and with respect to the initial data, 
we deduce that $A(\cdot)$ is continuous.

Moreover, we clearly have from \eqref{eq:Cauchypb}
$$
\tildeu \f{d\tildeu}{dy} = \tildec \f{d\tildec}{dy} \psi(-y).
$$
We have established the existence of $\tildec$ which is nonnegative and
nonincreasing on $(-1,0)$. Therefore, since $\psi(-y)\leq \psi(1)$
from assumption \eqref{hyppsi}, we deduce
$$
\f{d\tildeu^2}{dy} \geq \f{d\tildec^2}{dy} \psi(1).
$$
Integrating on $(-1,0)$ we obtain
$$
\tildeu^2(0)-\tildeu^2(-1) \geq \psi(1)(\tildec^2(0)-\tildec^2(-1))\geq -\psi(1).
$$
By definition $A(\sg)=\tildeu(-1)$ and with the boundary conditions in \eqref{eq:Cauchypb}, we conclude that the estimate $A(\sg)\leq \sqrt{\psi(1)}$ holds.
\qed

\subsection{Analytical example}

The Theorem~\ref{th:TWexist} can be illustrated thanks to an explicit traveling wave,  in the particular case
when 
\beq\label{eq:psicGc}
\psi(n)=\left\{\begin{array}{ll}
\lambda,&n=1,\\
\lambda n_c,&n<1,\\
0,&n=0.
\end{array}\right.
\qquad 
G(c)=\left\{\begin{array}{ll}
\ \ g_+,& \quad c>\bar{c},\\
-g_-,& \quad c<\bar{c}.
\end{array}\right.
\eeq
with $\lambda$, $n_c$, $g_+$, $g_-$ some positive constants and $n_c<1$. 
We can not apply directly Theorem \ref{th:TWexist} since the functions
$G$ and $\psi$ are not $C^1$. However, with this simple choice, we are 
able to manage explicit calculations and therefore state the following existence
result~:
\begin{theorem}\label{th:TWex}
Let $c_B>\bar{c}>0$ and let $\psi(n)$, $G(c)$ be chosen as in  \eqref{eq:psicGc}. 
There exists an unique monotone traveling wave for the Hele-Shaw model with 
nutrient i.e. an unique positive $\sigma$ and $R>0$ such that the system
\eqref{eq:TWg1}--\eqref{eq:jumprel} admits a nonnegative solution with $c' \geq 0$.
Moreover we have
\beq\label{exprsigma}
\sigma = R\;  \big(\sqrt{(g_++g_-) g_-}-g_-\big),
\eeq
where the value $R$ is obtained by solving \eqref{eqforR} below. The velocity $\sigma$ increases with respect to $c_B$ and decreases with respect to $\bar{c}$.
\end{theorem}
{\bf Proof.}
We use the same strategy as in the previous Section.
Let $\sg>0$ be given. We divide the whole real line into 
$I\cup II \cup III=(-\infty,0]\cup (0,R)\cup [R,+\infty)$.
\\
\\
$\bullet$ In $III=[R,+\infty)$, $n(x)=0$, $p(x)=0$, $c(x)=c_B$. 
\\[5pt]
$\bullet$ In $II=(0,R)$, $n(x)=1$ and $c(x)$ satisfies
$$
-\p_{xx}c+\lambda c=0, \qquad c(R)=c_B,\quad c'(R)=c'_R.
$$
As above, we have then that
$$
c(x)=c_B\cosh(\sqrt{\lambda}(x-R))+c'_R\sinh(\sqrt{\lambda}(x-R)).
$$
This function is positive and nondecreasing provided (see \eqref{boundcR})
\beq\label{eq:cR}
c_B \tanh(\sqrt{\lambda} R) \leq c'_R < \frac{c_B}{\tanh(\sqrt{\lambda}R)}.
\eeq
 
Consider the equation for $p$~: $-\p_{xx}p=G(c)$. On the one hand, we have $p(0)=0$ and 
$p'(0)= 0$ from \eqref{eq:jumprel}.
If we want $p$ to take nonnegative values, it has to be convex at $0$. 
Therefore from \eqref{eq:TWg4} we deduce $G(c(0))=-g_-<0$, 
i.e. $c(0)<\bar{c}$ from \eqref{eq:psicGc}.
Then we solve $-\p_{xx}p=-g_-$ with boundary conditions 
$p(0)=0$, $p'(0)=0$, which yields
$$
p(x)=\f{g_-}{2}x^2,\qquad x\in(0,x_1),\quad\mbox{where}\quad c(x_1)=\bar{c}.
$$
Since $c(x)$ is nondecreasing provided \eqref{eq:cR} is satisfied, 
$x_1>0$ and from $p(R)=0$, we have $R>x_1$. 
We can determine $x_1$ by $c(x_1)=\bar{c}$ such that
\beq\label{eq:x1}
c_B\cosh(\sqrt{\lambda}(x_1-R))+c'_R\sinh(\sqrt{\lambda}(x_1-R))
=\bar{c}.
\eeq
On the interval $x\in(x_1,R)$, $G(c(x))=g_+$ so that $-\p_{xx}p=g_+$.
Solving this equation with continuity of $p$ and  $p'$ at $x_1$,
we obtain
$$
p(x)=-\f{g_+}{2}(x-x_1)^2 +g_-x_1(x-x_1) + \f{x_1^2}{2}g_-,\qquad
x\in(x_1,R).
$$
Then, from $p(R)=0$, we have
\beq\label{eq:R}
-\f{g_+}{2}(R-x_1)^2 +g_-x_1(R-x_1) +\f{x_1^2}{2}g_-=0.
\eeq
Additionally, we recall the jump condition $\sigma=-p'(R)$ which gives
\beq\label{eq:sigma}
\sigma= g_+R-(g_++g_-)x_1.
\eeq
\\
$\bullet$ In the interval $I=(-\infty,0)$, $p(x)=0$. Since $c(x)\leq c(0)<\bar{c}$,
we have $-\sigma n'=-g_-n$. 
Together with the boundary conditions $n(0)=1$ (see \eqref{eq:jumprel}), we find
$$
n(x)= e^{\f{g_-}{\sigma}x}.
$$
Then $c(x)$ satisfies $-\p_{xx}c+\lambda n_cc=0$, together with the continuity of $c$
and $c'$ at $0$, we get
$$
c(x)=\f{1}{2}\Big(c_0+\f{c'_0}{\xi}\Big)e^{\xi x}+
\f 12 \Big(c_0-\f{c'_0}{\xi}\Big)e^{-\xi x}, 
\qquad \mbox{ with } \xi=\sqrt{\lambda n_c},
$$
where the expressions of $c_0$ and $c'_0$ are given in \eqref{eq:c0}--\eqref{eq:cprime0}.
In order to have $c$ bounded on $(-\infty,0)$, the coefficient in front of $e^{-\xi x}$
has to be zero~:
\begin{equation}
  \label{eq:coe0}
  c_0 - \f{c'_0}{\xi}=0.
\end{equation}
Then we have $\dis c(x) = c_0 e^{\xi x}$.
We observe that $c(x)\geq 0$ for $x\in (-\infty,0)$
provided \eqref{eq:cR} is satisfied.
\\

With this construction, we have to solve the three equations  \eqref{eq:x1}, \eqref{eq:R} and
\eqref{eq:coe0}, for the  three unknowns~:
$c'_R$, $x_1$, $R$. In fact, we notice that all these equations are independent
of $\sigma$; more precisely, compared to the proof of Theorem \ref{th:TWexist},
we have $A(\sg)=\xi$ from \eqref{eq:coe0}. 
Then \eqref{eq:sigma} will give the velocity $\sigma$.
First, we recover from \eqref{eq:coe0},
\beq\label{eq:cL1}
c'_R = c_B \f{\xi\cosh(\sqrt{\lambda}R)+\sqrt{\lambda}\sinh(\sqrt{\lambda} R)}
{\xi\sinh(\sqrt{\lambda} R)+\sqrt{\lambda}\cosh(\sqrt{\lambda}R)}.
\eeq
This identity is equivalent to \eqref{eq:c'R} with $A(\sigma)=\xi$.
We notice that this expression implies that \eqref{eq:cR} is satisfied.
Then we solve equation \eqref{eq:R} and obtain
$$
x_1 = R \Big(1-\sqrt{\f{g_-}{g_++g_-}} \Big).
$$
We set
$$
\alpha := \sqrt{\f{g_-}{g_++g_-}} \in (0,1).
$$
Finally, substituting $x_1=(1-\alpha) R$ and \eqref{eq:cL1} into \eqref{eq:x1},
we obtain the following nonlinear equation for $R>0$~:
\beq\label{eqforR}
{\cal F}(R):=\f{\bar{c}}{c_B}\Big(\cosh(\sqrt{\lambda} R) + \sqrt{n_c} \sinh(\sqrt{\lambda} R)\Big) 
-\sqrt{n_c} \sinh(\sqrt{\lambda}(1-\alpha) R) -\cosh(\sqrt{\lambda}(1-\alpha) R) = 0.
\eeq
The function ${\cal F}$ is $C^\infty(\R_+)$ and we prove below that there exists
an unique positive root for this function. Then we obtain the 
expression \eqref{exprsigma} for the velocity from equation \eqref{eq:sigma}.

Let us state the existence of an unique positive root for ${\cal F}$.
We have ${\cal F}(0)=(\f{\bar{c}}{c_B}-1)< 0$, and 
${\cal F}(R) \to +\infty$ when $R\to +\infty$.
Then there exists a positive root for ${\cal F}$.
To prove the uniqueness of this root, let us first choose $k\in\N^*$
such that $(1-\alpha)^k\leq \f{\bar{c}}{c_B}$ and $(1-\alpha)^{k-1}>\f{\bar{c}}{c_B}$.
For every $i\in \N$, we compute the $i$th derivative of ${\cal F}$ from 
\eqref{eqforR}~:
$$
\begin{array}{ll}
\dis {\cal F}^{(2i)}(R)= 
&\dis \sqrt{\lambda}^{{2}i}\Big(\f{\bar{c}}{c_B}\big(\cosh(\sqrt{\lambda} R) + \sqrt{n_c} \sinh(\sqrt{\lambda} R)\big) \\[2mm]
&\dis \quad -(1-\alpha)^{2i}\big(\sqrt{n_c} \sinh(\sqrt{\lambda}(1-\alpha) R) +\cosh(\sqrt{\lambda}(1-\alpha) R)\big)\Big).\\
\dis {\cal F}^{(2i+1)}(R)= 
&\dis \sqrt{\lambda}^{2i+1}\Big(\f{\bar{c}}{c_B}\big(\sinh(\sqrt{\lambda} R) + \sqrt{n_c} \cosh(\sqrt{\lambda} R)\big) \\[2mm]
&\dis \quad -(1-\alpha)^{2i+1}\big(\sqrt{n_c} \cosh(\sqrt{\lambda}(1-\alpha) R) +\sinh(\sqrt{\lambda}(1-\alpha) R)\big)\Big).
\end{array}
$$
For all $i\geq 0$, we have ${\cal F}^{(i)}(R)\to +\infty$ as $R\to +\infty$.
Thanks to our choice for $k$, we have ${\cal F}^{(k)}(R)> 0$ for all $R>0$
and ${\cal F}^{(k-1)}(0)<0$. 
Thus ${\cal F}^{(k-1)}$ is increasing and admits a unique root denoted $r_{k-1}$.
Therefore, if $k=1$ we have proved the result. 
Else, we have then ${\cal F}^{(k-2)}$ is a nonincreasing function on $(0,r_{k-1})$
and an increasing function on $(r_{k-1},+\infty)$.
We have from our choice of $k$ that ${\cal F}^{(k-2)}(0)<0$, then ${\cal F}^{(k-2)}(r_{k-1})<0$. 
Since ${\cal F}^{(k-2)}(R)\to +\infty$ as $R\to +\infty$, we conclude that
${\cal F}^{(k-2)}$ admits an unique positive root $r_{k-2}$.
Therefore, if $k=2$ we have proved the result, else we continue
and by induction we prove by the same token that 
${\cal F}$ is nonincreasing on $(0,r_1)$, increasing on
$(r_1,+\infty)$ with ${\cal F}(r_1)\leq {\cal F}(0)<0$ and
$\lim_{R\to +\infty} {\cal F}(R)=+\infty$.
Thus ${\cal F}$ admits an unique positive root.

Moreover, from the expression of ${\cal F}$ in \eqref{eqforR}, we deduce clearly
that $\dis \f{\p {\cal F}}{\p c_B} < 0$ and $\dis \f{\p {\cal F}}{\p \bar{c}} > 0$.
Thus we deduce that $R$ increases with respect to $c_B$ and decreases with 
respect to $\bar{c}$. From \eqref{exprsigma}, we conclude the proof.
\qed

\bigskip

The situation here is simpler than in Theorem \ref{th:TWexist}. Indeed, we can compute $R$ as the root of an explicit function  \eqref{eqforR} which only depends on the parameters. Then we have an analytic 
expression both of the shape of the wave and of the velocity. 

%
%

Some consequences of Theorem \ref{th:TWex} are compatible with the biological intuition. For example, the traveling velocity  increases with the nutrient density outside of the tumor, and decreases with
$\bar{c}$ which is the critical nutrient density for the growth of tumor cells.

\section{Traveling waves for the {\it in vivo} case in one dimension}
\label{sec:invivo}

The only difference between {\it in vivo} and {\it in vitro} models, is the equation for the nutrient availability.
The one dimensional {\it in vivo} model reads~:
\beq\label{eq:vitro}
\left\{
\begin{array}{ll}
-\sigma \p_x n =n G(c), 
& \mbox{in }   \mathbb{R}\setminus(0,R),
\\[1mm]
n=1&\mbox{in }(0,R),
\\[1mm]
-\p_{xx}c+\psi(n)c=\ind{n=0}(c_B-c),& x\in \R,
\\[1mm]
\dis \lim_{x\to + \infty} c(x)=c_B,
 \\[1mm]
-\p_{xx}p=G(c),& \mbox{in } (0,R),
\\[1mm]
p(0)= p(R)=0.
\end{array}
\right.
\eeq
This system is still complemented with the jump condition \eqref{eq:jumprel}.
Then, we can adapt the result in Section~\ref{sec:TW} and 
prove the following theorem.

\begin{theorem}\label{th:TWvitro}
Let $\dis 0<\bar{c}<\frac{c_B}{2}$ be given.
Assume that $G$ and $\psi$ are $C^1$ function satisfying \eqref{hypG}--\eqref{hyppsi}.
Then, there exist $\sigma>0$ and $R>0$ such that the system 
\eqref{eq:vitro} together with the jump condition \eqref{eq:jumprel} admits a solution with $c$ 
increasing on $\mathbb{R}$; moreover
$n$ is increasing on $(-\infty,0]$ and $\lim_{x\to -\infty} n(x) = 0$.
\end{theorem}
{\bf Proof.}
The idea of the proof is the same as the proof of Theorem \ref{th:TWexist}.
We will only give the main changes and do not provide all the details of the 
calculations.
\\
{\em Step 1. The piecewise construction of the wave.} 
Let us assume that $\sigma >0$ is given. We construct the solution for $R>0$ 
following the proof of Theorem \ref{th:TWexist}. 
We just explain below the main changes~:
\begin{itemize}
\item On $III=(R,+\infty)$, the only change is for $c$~:
$c(x)=c_B -c'_Re^{\sqrt{\psi(1)}(R-x)}$.
\item On $II=[0,R]$, the equations are the same but the boundary conditions
for $c$ at $x=R$ changes. Then we get
\beq\label{eq:cIIvitro}
c(x)=c_B \cosh\big(\sqrt{\psi(1)}(x-R)\big)-c'_R e^{-\sqrt{\psi(1)}(x-R)},
\eeq
which is nondecreasing and positive provided
\beq\label{boundcRvitro}
\frac{c_B}{2}\Big(1-e^{-2\sqrt{\psi(1)}R}\Big) \leq c'_R < \frac{c_B}{2} 
\Big(1+e^{-2\sqrt{\psi(1)}R}\Big).
\eeq
The pressure is given by expression \eqref{pII}; we deduce, from the boundary condition 
$p(R)=0$, a similar condition for $R$ as in \eqref{eq1}
\beq\label{eq1vitro}
\int_0^1 s G\Big(c_B \cosh\big(\sqrt{\psi(1)}Rs\big)-c'_R e^{\sqrt{\psi(1)}Rs}\Big)\,ds=0.
\eeq
\item On $I=(-\infty,R)$, applying Lemma \ref{lem1}, 
in order to have the existence of a solution of the system on $I=(-\infty,0)$ such that $c$ is nonnegative
and nondecreasing,
we introduce a nonnegative
function $A(\sg)$ such that $c'(0)=A(\sg)c(0)$.
Injecting the expressions of $c(0)$ and $c'(0)$ from \eqref{eq:cIIvitro} into this identity, we deduce
\beq\label{exprc'R}
c'_R= \frac{c_B}{2}\Big(1+\frac{A(\sg)-\sqrt{\psi(1)}}{A(\sg)+\sqrt{\psi(1)}} 
e^{-2\sqrt{\psi(1)}R}\Big).
\eeq
We deduce then that \eqref{boundcRvitro} is satisfied.
\end{itemize}
\smallskip
{\em Step 2. The nonlinear equation for $R$ and $c'_R$.}
To conclude the construction, it suffices to show that the nonlinear equation
\eqref{eq1vitro} admits a solution $R_\sigma$. Let us first introduce the notation
$$
\gamma_\sg(R,s) := c_B \cosh(\sqrt{\psi(1)}Rs)-c'_R e^{\sqrt{\psi(1)}Rs}, \qquad \mbox{for }s\in(0,1).
$$
With \eqref{exprc'R}, we deduce that
\beq\label{eq:gammaRs}
\gamma_\sg(R,s) = \f{c_B}{2} \Big( e^{-\sqrt{\psi(1)}Rs} + \f{\sqrt{\psi(1)}-A(\sg)}{A(\sg)+\sqrt{\psi(1)}} e^{\sqrt{\psi(1)}R(s-2)}\Big).
\eeq
Then \eqref{eq1vitro} rewrites $\int_0^1 sG(\gamma_\sg(R,s))ds=0$.
From the estimate $A(\sg)\leq \sqrt{\psi(1)}$ in Lemma \ref{lem1}, we deduce straightforwardly that $R\mapsto \gamma_\sg(R,s)$ is decreasing
for $s\in(0,1)$.
With \eqref{eq:gammaRs} we have $\dis \gamma_\sg(0,s)=c_B\f{\sqrt{\psi(1)}}{A(\sg)+\sqrt{\psi(1)}}$
and $\dis \lim_{R\to + \infty} \gamma_\sg(R,s)=0$.
Since $G(0)<0$, we deduce that $\int_0^1 sG(\gamma_\sg(R,s))ds$ takes negative value for
large $R$. 
From Lemma \ref{lem1}, we have $A(\sg) \leq \sqrt{\psi(1)}$.
Then, $\gamma_\sg(0,s)\geq c_B/2$.
By assumption $\bar{c}< c_B/2$, we deduce then from \eqref{hypG} that 
$\int_0^1 sG(\gamma_\sg(0))ds>0$. 
By continuity and monotonicity, there exists a unique $R_\sg >0$ such that
\eqref{eq1vitro} is satisfied.
Moreover, as above, the mapping $\sg\mapsto R_\sg$ is continuous.

Then for all $\sigma>0$, we have constructed $R_\sg$ and a solution $(n_\sg,c_\sg,p_\sg)$
to the system \eqref{eq:vitro}. 
\\
\\
{\em Step 3. The nonlinear equation for $\sg$.} 
It suffices now to prove that there exists $\sg>0$ 
such that $\sigma = -p'_\sg(R_\sg)$, which can be rewritten~:
\beq\label{eqsgvitro}
\sigma = R_\sg \int_0^1 G\big(c_B\cosh(\sqrt{\psi(1)}R_\sg s)-c'_R e^{\sqrt{\psi(1)}R_\sg s}\big)\,ds.
\eeq
Moreover, since we have for all $R>0$,
$$
\gamma_\sg(R)\leq \f{c_B}{2} e^{-\sqrt{\psi(1)}Rs} \Big(1+e^{-\sqrt{\psi(1)}Rs}\Big),
$$
we deduce that $R_\sg\leq R_b$, where $R_b$ is the unique solution of the equation
$$
\int_0^1 s G\Big(\f{c_B}{2} e^{-\sqrt{\psi(1)}Rs} \big(1+e^{-\sqrt{\psi(1)}Rs}\big)\Big)\,ds=0.
$$
By the same token as in the proof of Theorem \ref{th:TWexist}, we have that
the right hand side of \eqref{eqsgvitro} is bounded and admits a positive
and finite limit when $\sigma\to 0$. Therefore, 
by continuity, there exists $\sg>0$ such that \eqref{eqsgvitro} is satisfied.
\qed

\begin{remark}
Notice that the condition $\bar{c}<c_B/2$ can be understood from the requirement that 
$c(0)<\bar{c}<c(R)$. From the expression of $\gamma_\sg(R,s)$ in \eqref{eq:gammaRs}, 
$$
c(0)=c_Be^{-\sqrt{\psi(1)}R}\f{\sqrt{\psi(1)}}{A(\sigma)+\sqrt{\psi}},\qquad 
c(R)=\f{c_B}{2}\left(1-\f{A(\sigma)-\sqrt{\psi(1)}}{A(\sigma)+\sqrt{\psi(1)}}e^{-2\sqrt{\psi(1)}R}\right).
$$
For fixed $A(\sigma)$, both $c(0)$ and $c(R)$ decrease with $R$. Moreover, 
when $R\to+\infty$, $c(0)\to 0$ and $c(R)\to c_B/2$, while when $R\to 0$, $c(0)=c(R)=c_B\f{\sqrt{\psi(1)}}{A(\sigma)+\sqrt{\psi(1)}}$.
Therefore, when $\bar{c}<c_B/2$, $\bar{c}<c(R)$ holds for any $R$. Furthermore, $R$ big enough gives $c(0)<\bar{c}$.

This condition confirms our numerical observations in Section
\ref{sec:num}. In fact, for the {\it in vivo} model,  we have noticed that
traveling waves do not exist when we use the expression in \eqref{eq:ex1parameter} 
for the growth function for which the condition $\bar{c}<c_B/2$ is violated. 
But if we choose \eqref{num:Ginvivo} instead, numerical simulations
exhibits traveling waves (see Fig.~\ref{fig:ex1_TWv}).

\end{remark}
\bigskip
{\bf Analytical example.}

Following the previous section, we consider the case where the nutrient consumption
$\psi$ and the growth term $G$ are given by \eqref{eq:psicGc}. Then, we have the
\begin{proposition}
Let $\dis 0<\bar{c}<\f{c_B}{2}$ and assume that \eqref{eq:psicGc} holds true
with $n_c<1$.
Then there exists a unique traveling wave.
Moreover, its velocity is given by
$$
\sigma = R\big(\sqrt{(g_++g_-)g_-}-g_-\big),
$$
where $R$ is the solution to \eqref{eqforRvitro}. Moreover, $\sigma$ is nondecreasing 
with respect to $c_B$ and nonincreasing with respect to $\bar{c}$.
\label{prop:analytic_vivo}
\end{proposition}
{\bf Proof.}
As in Theorem \ref{th:TWvitro}, in the region $[0,R]$, the nutrient concentration $c(x)$
can be given by \eqref{eq:cIIvitro} 
with $\psi(1)=\lambda$.
As has been stated in the proof of Theorem \ref{th:TWex} that when $\psi(n)$, $G(n)$
are as in \eqref{eq:psicGc}, we have $A(\sg)=\sqrt{\lambda n_c}$.
Thus identity \eqref{exprc'R} becomes
$$
c'_R = \frac{c_B}{2}\Big(1-\f{1-\sqrt{n_c}}{1+\sqrt{n_c}}
e^{-2\sqrt{\lambda} R}\Big).
$$

From the proof of Theorem \ref{th:TWex}, 
let $c(x_1)=\bar{c}$, the equation for the pressure indicates that 
$x_1=(1-\sqrt{\f{g_-}{g_++g_-}})R=(1-\alpha)R$ with $\alpha = \sqrt{\f{g_-}{g_++g_-}}$.
Then, from \eqref{eq:cIIvitro}, the nonlinear equation for the unknonw $R$ is~:
$$
c_B \cosh(\sqrt{\lambda} \alpha R) - \frac{c_B}{2}\Big(1-
\f{1-\sqrt{n_c}}{1+\sqrt{n_c}}e^{-2\sqrt{\lambda} R}\Big)
e^{\alpha\sqrt{\lambda} R} = \bar{c}, \qquad \alpha = \sqrt{\f{g_-}{g_++g_-}}.
$$
After simplifications, this latter equation rewrites
\beq\label{eqforRvitro}
{\cal F}(R):= \frac{c_B}{2}\Big(e^{-\alpha\sqrt{\lambda}R}
+\f{1-\sqrt{n_c}}{1+\sqrt{n_c}}e^{(\alpha-2)\sqrt{\lambda} R}\Big)
= \bar{c}, \qquad \alpha = \sqrt{\f{g_-}{g_++g_-}}.
\eeq
We have $\dis {\cal F}(0)=\f{c_B}{1+\sqrt{n_c}} > \f{c_B}{2} > \bar{c}$ and
$\dis \lim_{R\to +\infty}{\cal F}(R)=0$. 
Moreover, ${\cal F}$ is nonincreasing, then there exists a unique 
solution $R>0$ of \eqref{eqforRvitro} which is independent of $\sg$.
Finally, the velocity $\sigma$ is given by \eqref{eq:sigma}.
Furthermore, we clearly have that $c_B\mapsto {\cal F}(R)$ is a nondecreasing function,
which gives monotonicity of $R$ with respect to $c_B$ and then of $\sigma$ with respect
to $c_B$ and to $\bar{c}$,  as stated in Theorem \ref{prop:analytic_vivo}.
\qed

\begin{remark}[Comparison of the velocity]
We notice that, for this analytical example, in both {\it in vitro} and 
{\it in vivo} cases, the expression of the velocity $\sigma$ with respect 
to $R$ is the same.
The difference being that the radius of the tumor $R_{vitro}$ and $R_{vivo}$ 
do not satisfy the same nonlinear equation.
However we can compare them by noticing that we can rewrite \eqref{eqforR} as
$$
\bar{c} = c_B \f{\dis e^{-\alpha\sqrt{\lambda} R_{vitro}}+\f{1-\sqrt{n_c}}{1+\sqrt{n_c}}e^{(\alpha-2)\sqrt{\lambda}R_{vitro}}}
{\dis 1+\f{1-\sqrt{n_c}}{1+\sqrt{n_c}}e^{-2\sqrt{\lambda}R_{vitro}}}
> \f{c_B}{2}\Big( e^{-\alpha\sqrt{\lambda} R_{vitro}}+\f{1-\sqrt{n_c}}{1+\sqrt{n_c}}e^{(\alpha-2)\sqrt{\lambda}R_{vitro}}\Big).
$$
We recognize the expression in \eqref{eqforRvitro} in the right hand side.
Thus taking the notation in \eqref{eqforRvitro}, we have 
${\cal F}(R_{vitro})<\bar{c}={\cal F}(R_{vivo})$.
Since ${\cal F}$ is nonincreasing, we deduce that $R_{vivo}\leq R_{vitro}$.
Therefore, we can conclude that $\sg_{vivo}\leq \sg_{vitro}$.
\end{remark}

\section{Conclusions and perspectives}

Numerical solutions of models of tumor growth, based on a free boundary problem coupled to a diffusion equation for the nutrient, exhibit complex structures in addition to the usually observed proliferative rim and necrotic core. Remarkable are the possible radial instability in two dimensions and the profiles for the cell number density and the local pressure.  We have been able to explain the latter with a study of  the traveling waves for the model under consideration. The structure of the tumor is then described thanks to (nearly) analytical formulas:  a sharp front separates the healthy tissue from the proliferative rim  where the density
is close to his highest value;  the proliferative rim undergoes  a smooth transition (in particular the pressure is differentiable) to the necrotic core where the mechanical pressure vanishes. 
\\

Our calulations are tractable because we used particular nonlinearities; a general study of the existence of traveling waves is still to be done for general nonlinearities and for finite values of $\gamma$ in equation \eqref{eq:pqn} with $p(n)=n^\gamma$. 
\\

We have  not studied the radial instability appearing in the 
numerical simulations, a complex phenomena already observed on a simpler system in \cite{Kessler_Levine}. In order to analyze this instability, possible strategies are those developed in \cite{BenAmar} or \cite{KPV}. We will come back on this question in a forthcoming work. 

\bigskip

\noindent
{\bf Acknowledgements:} The authors would like to thank partial support
from Campus France program~: Xu GuangQi Hubert Curien program n$^o 30043VM$
{\it PDE models for cell self-organization}. MT is supported by Shanghai Pujiang Program 13PJ1404700 and NSFC11301336. BP and NV are also supporter by ANR KIBORD, No ANR-13-BS01- 0004-01.
%
%

\end{document}